\documentclass[11pt,letterpaper]{amsart}
\pdfoutput=1
\usepackage{arxiv_math}
\usepackage{extarrows,multirow,colortbl,hhline,transparent}
\usepackage{amssymb,amsmath,array,diagbox}
\usepackage{xypic}
\usepackage{import}

\author[G.\ Guth]{Gary Guth}
\address{Department of Mathematics, Stanford University, Stanford, CA 94305 USA}
\email{gmguth@stanford.edu}

\author[K.\ Hayden]{Kyle Hayden}
\address{Department of Mathematics and Computer Science, Rutgers University, Newark, NJ 07102, USA}
\email{kyle.hayden@rutgers.edu}

\author[S.\ Kang]{Sungkyung Kang}
\address{Mathematical Institute, University of Oxford, Andrew Wiles Building, Radcliffe Observatory Quarters, Woodstock Road, Oxford, OX2 6GG, United Kingdom}
\email{sungkyung3838@gmail.com}

\author[J.\ Park]{JungHwan Park}
\address{Department of Mathematical Sciences, KAIST, Daejeon 34141, Republic of Korea}
\email{jungpark0817@kaist.ac.kr}

\definecolor{cellgray}{gray}{0.875}

\begin{document}
\title{Doubled Disks and Satellite Surfaces}
\begin{abstract}
    Conjecturally, a knot is slice if and only if its positive Whitehead double is slice. We consider an analogue of this conjecture for slice disks in the four-ball: two slice disks of a knot are smoothly isotopic if and only if their positive Whitehead doubles are smoothly isotopic. We provide evidence for this conjecture, using a range of techniques. More generally, we consider when isotopy obstructions persist under satellite operations. In particular, we show that obstructions coming from knot Floer homology, Seiberg-Witten theory, and Khovanov homology often behave well under satellite operations. 

    %More generally, we consider when isotopy obstructions persist under satellite operations. We show that obstructions coming from knot Floer homology, Seiberg-Witten theory, and Khovanov homology often behave well under satellite operations. We apply this strategy to give a systematic method for constructing vast numbers of exotic disks in the four-ball (which can be upgraded to stably-exotic disks using the same techniques) and to produce an infinite family of pairwise exotic slice disks with boundary $\Wh(4_1 \# 4_1)$.
     We apply these strategies to give a systematic method for constructing vast numbers of exotic disks in the four-ball, including the first infinite family of pairwise exotic slice disks. These same techniques are then upgraded to produce exotic disks that remain exotic after  any prescribed number of internal stabilizations. Finally, we show that the branched double covers of certain stably-exotic disks become diffeomorphic after a single stabilization with $S^2 \times S^2$, hence stabilizing them yields exotic surfaces that have diffeomorphic branched covers.
\end{abstract}
\thanks{GG is partially supported by NSF Grant DMS-2204214 and a Simons Collaboration Grant on New Structures in Low Dimensional Topology. KH is supported by NSF grant DMS-2243128. JP is partially supported by Samsung Science and Technology Foundation (SSTF-BA2102-02) and the POSCO TJ Park Science Fellowship.}

\maketitle

%\tableofcontents 

\section{Introduction}\label{sec: intro}
%\footnote{\JP{Feel free to change the title!}\GG{I'm temporarily putting in an introduction for my letter writers. We can certainly replace any of this later!}}

Given a knot $K$ in the 3-sphere and a pattern knot $P$ embedded in the solid torus, the \emph{satellite of $K$ with pattern $P$} is defined by cutting out a neighborhood of $K$ and filling it with the solid torus containing $P$. We will denote the resulting knot $P(K)$. 

\begin{figure}\center
\def\svgwidth{.8\linewidth}%% Creator: Inkscape 1.2 (dc2aeda, 2022-05-15), www.inkscape.org
%% Accompanies image file '946-doubles.pdf' (pdf, eps, ps)
\begingroup%
  \makeatletter%
  \providecommand\color[2][]{%
    \errmessage{(Inkscape) Color is used for the text in Inkscape, but the package 'color.sty' is not loaded}%
    \renewcommand\color[2][]{}%
  }%
  \providecommand\transparent[1]{%
    \errmessage{(Inkscape) Transparency is used (non-zero) for the text in Inkscape, but the package 'transparent.sty' is not loaded}%
    \renewcommand\transparent[1]{}%
  }%
  \providecommand\rotatebox[2]{#2}%
  \newcommand*\fsize{\dimexpr\f@size pt\relax}%
  \newcommand*\lineheight[1]{\fontsize{\fsize}{#1\fsize}\selectfont}%
  \ifx\svgwidth\undefined%
    \setlength{\unitlength}{830.19866703bp}%
    \ifx\svgscale\undefined%
      \relax%
    \else%
      \setlength{\unitlength}{\unitlength * \real{\svgscale}}%
    \fi%
  \else%
    \setlength{\unitlength}{\svgwidth}%
  \fi%
  \global\let\svgwidth\undefined%
  \global\let\svgscale\undefined%
  \makeatother%
  \begin{picture}(1,0.79941099)%
    \lineheight{1}%
    \setlength\tabcolsep{0pt}%
    \put(0,0){\includegraphics[width=\unitlength,page=1]{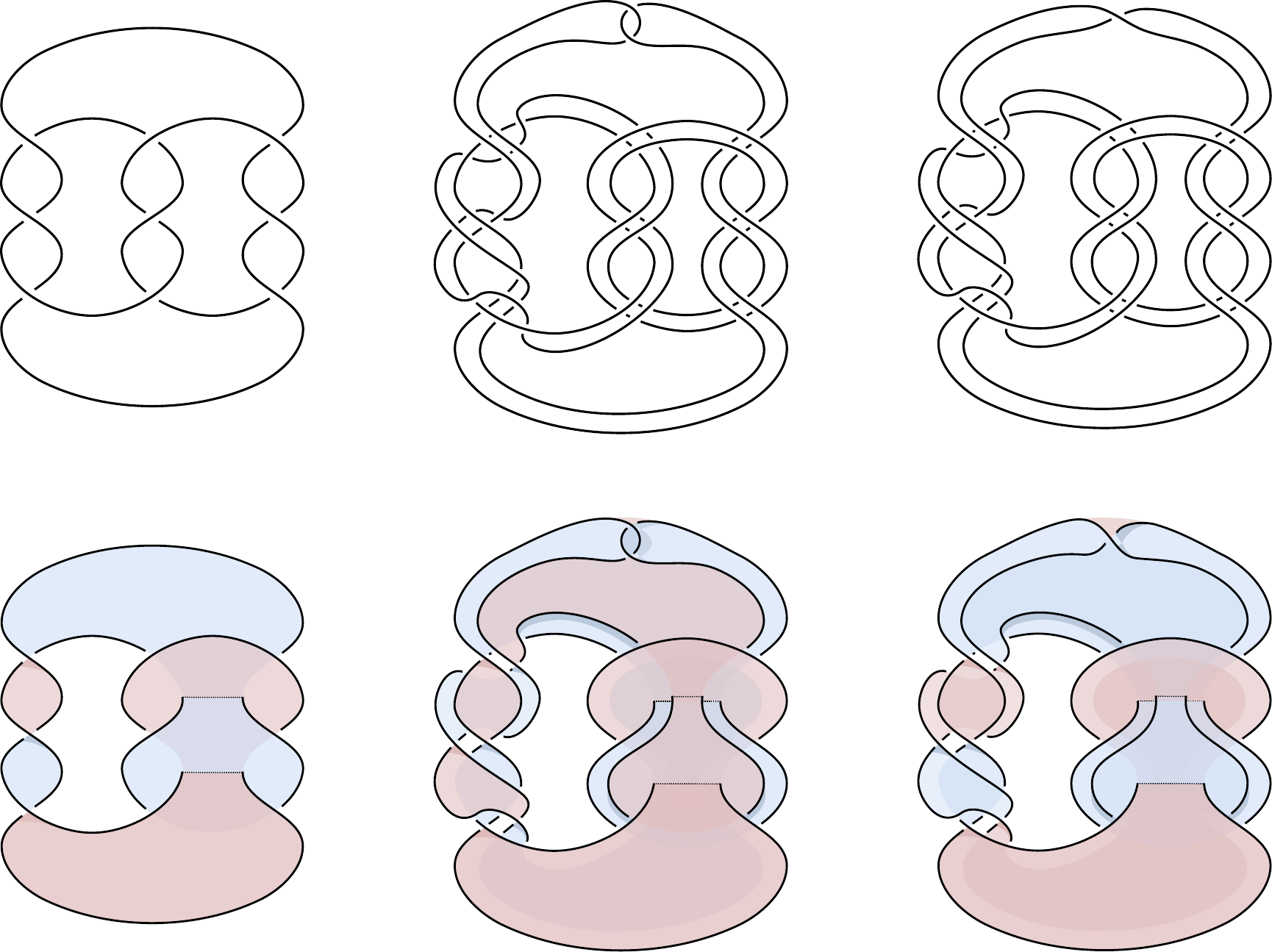}}%
    \put(0.10037481,0.03022288){\color[rgb]{0,0,0}\makebox(0,0)[lt]{\smash{\begin{tabular}[t]{l}{\small $D$}\end{tabular}}}}%
    \put(0.44307477,0.00823589){\color[rgb]{0,0,0}\makebox(0,0)[lt]{\smash{\begin{tabular}[t]{l}{\small $\operatorname{Wh}(D)$}\end{tabular}}}}%
    \put(0.84960459,0.00775949){\color[rgb]{0,0,0}\makebox(0,0)[lt]{\smash{\begin{tabular}[t]{l}{\small $D_{2,1}$}\end{tabular}}}}%
    \put(0,0){\includegraphics[width=\unitlength,page=2]{946-doubles.pdf}}%
  \end{picture}%
\endgroup%

\caption{A slice disk for $m(9_{46})$, together with its Whitehead double and $(2,1)$-cable.}\label{fig:946-doubles}
\end{figure}

Satellite operations are ubiquitous in the study of knots in dimension three, but also have natural four dimensional extensions. Given a concordance $C$ from $K_0$ to $K_1$, there is satellite concordance from $P(K_0)$ to $P(K_1)$, denoted $P(C)$, which is obtained by cutting out a neighborhood of the concordance and gluing back in the pattern cylinder $(S^1\times D^1\times [0, 1], P \times [0, 1])$. When $K_0$ is the unknot and the pattern $P$ is unknotted, the satellite construction uniquely determines a satellite slice disk for $P(K_1)$; see \Cref{fig:946-doubles} for examples. These ideas extend further to higher-genus surfaces and surfaces in other 4-manifolds (e.g., \cite{shinohara,hkkmps,Hayden_Kim_Miller_Park_Sundberg_Kh_exotic_seifert}). %Informally, one fixes a tubular neighborhood $N$ of a surface $S$ in a 4-manifold and replaces the base $S$ with another surface in $N$. In this paper, we are chiefly concerned with the case of surfaces in $B^4$, where the tubular neighborhood is given by $N\cong S \times D^2$ , and the satellite surface is typically constructed from disjoint sections of the form $S \times \{pt\}$ joined by bands attached to their boundaries in $\partial S \times D^2$ in $S^3$.

The question of how knot invariants behave under satellite operations is well-studied. We aim to further the analogous study of the behavior of surface invariants under satellite operations, and tease out some of the similarities in the two contexts. We conjecture that many satellite operations are ``injective''; we show, that for large families of satellite patterns, disks which can be smoothly distinguished often have satellite disks which are also smoothly distinguishable. Matters are quite different in the topological category: by the work of Conway and Powell \cite{conway2021embedded}, under suitable satellite operations, disks which were initially topologically distinct can become topologically isotopic. This stark contrast leads us to consider the ways in which satellite operations can be used in the study of exotic phenomena.

%\KH{[just sticking this here for now while being too lazy to integrate it:]} 
In particular, we see that satellite operations can help address a fundamental dilemma:  finding constructions that are mild enough to preserve the topological type of an object yet are significant enough to cause a measurable change in its smooth type. %(And, in practice, we must detect this change using a computable smooth invariant.) 
For example, Fintushel and Stern's \emph{rim surgery} operation \cite{fintushel1997surfaces} aims to strike this balance directly, and variations of this operation have led to virtually all known examples of infinite families of exotic surfaces. An alternative is to begin with a violent operation that changes both the smooth \emph{and} topological type, and then apply a second operation (e.g., Whitehead doubling) to erase the topological difference between the surfaces.

In this paper, we study slice disks that are \emph{exotically knotted} (i.e., topologically isotopic rel boundary but not smoothly isotopic rel boundary)  using three different techniques: knot Floer homology, Khovanov homology, and Seiberg-Witten invariants. We briefly outline their strengths and weaknesses here. Knot Floer homology supports the most broadly applicable techniques among the three, and indeed it produces the most general statements in this paper thanks to tools coming from bordered Floer homology \cite{LOT_bordered_HF}. In contrast, making general statements about surface invariants in Khovanov homology is quite difficult, especially given its lack of analogous bordered techniques.  However, it has the advantage of being inherently combinatorial, which allows us to make simple computations. %Interestingly, Khovanov homology can also treat surfaces of non-minimal genus, which is more difficult using knot Floer techniques for grading reasons. 
Finally, we make essential use of the fact that the Seiberg-Witten invariants are defined over $\mathbb{Z}$ coefficients; this allows us to distinguish infinite families of pairwise exotic disks bounded by the same knot.

%In this paper, we will use three different techniques to detect exotic disks: knot Floer homology, Khovanov homology, and Seiberg-Witten invariants. We briefly outline their strengths and weaknesses here. Knot Floer homology provides the most broadly applicable techniques among the three, and indeed it produces the most general statements in this paper thanks to tools coming from bordered Floer homology \cite{LOT_bordered_HF}. In contrast, making general statements about surface invariants in Khovanov homology is quite difficult without cut-and-paste type arguments. However, it has the advantage of being inherently combinatorial, which allows us to make simple computations. Interestingly, Khovanov homology can also treat surfaces of non-minimal genus, which is almost impossible using knot Floer techniques. Finally, we produce an infinite family of pairwise exotic disks with the same boundary which we distinguish using Seiberg-Witten invariants, which have the advantage of being definable over $\mathbb{Z}$ coefficients.

%\KH{[Perhaps say something about how each obstruction we consider exhibits a strength missing from the other obstructions: knot Floer homology allows us to establish broadly-applicable sufficient conditions for producing satellite disks that are exotic (and even stably exotic); Khovanov homology allows us to make simple computations and treat surfaces of arbitrary genus; Seiberg-Witten invariants allow us to distinguish infinite families of pairwise exotic disks.]}

\subsection{Satellite Operators}\label{sec:injective satellite}

By the work of Freedman \cite{freedman,freedman_quinn_top_4mflds}, any knot in $S^3$ with trivial Alexander polynomial knot is topologically slice.  The Alexander polynomial behaves well with respect to satellite operations: if $P(K)$ is a satellite knot with pattern $P$ and companion $K$, the Alexander polynomial of $P(K)$ can be computed as
\[
\Delta_{P(K)}(t) = \Delta_P(t)\cdot \Delta_K(t^w),
\]
where the pattern $P$ represents $w$ times the generator for $H_1(S^1\times D^2)$. As a consequence, the positive Whitehead double of a knot $K$ is always topologically slice. %(though, in many instances, not smoothly slice). An open conjecture states that the sliceness of a knot $K$ is directly related to the sliceness of its positive Whitehead double.
This contrasts with the smooth setting, where the following remains a major open problem: 

\begin{conjecture}[{\cite[Problem~1.38]{kirby_problem_list}}]\label{conj: whitehead sliceness}
A knot $K$ is smoothly slice if and only if its positive Whitehead double is smoothly slice. 
\end{conjecture}

In \cite{hedden_whitehead_doubles}, Hedden provided strong evidence for this conjecture by showing that if $\tau(K) > 0$ then $\tau(\Wh(K))$ is nonzero. In  particular, if knot Floer homology obstructs $K$ from being smoothly slice it will obstruct $\Wh(K)$ from being smoothly slice as well. 

Inspired by the situation in three dimensions, we consider the analogous case of slice disks in the four-ball. Conway and Powell \cite{conway_powell_characterisation_homotopy_ribbon} showed that if $D_1$ and $D_2$ are slice disks for a knot $K$, then they are topologically isotopic if 
\[
\pi_1(B^4\setminus \nu D_1) \cong \Z \cong \pi_1(B^4\setminus \nu D_2).
\]
As a consequence, the positive Whitehead doubles of any two slice disks for a knot $K$ are topologically isotopic (see Section \ref{sec: background}). Given Conjecture \ref{conj: whitehead sliceness}, it is natural to ask whether the smooth isotopy class of a slice disk is determined by that of its Whitehead double.

\begin{conjecture}\label{conj: whitehead doubles for disks}
Two slice disks for a knot $K$ are smoothly isotopic if and only if their positive Whitehead doubles are smoothly isotopic. 
\end{conjecture}

%An unknotted pattern $P$ induces a map
%\[
%\{\text{isotopy classes of slice disks for } K\} \xra{P} \{\text{isotopy classes of slice disks for } P(K)\} 
%\]
%in the obvious way, by taking a slice disk for $K$ to the satellite disk for $P(K)$ determined by the pattern. In this language, we can reframe Conjecture \ref{conj: whitehead doubles for disks} as asking whether the positive Whitehead double operator (which we will denote $\Wh$) is injective. We provide evidence for this conjecture (as well as evidence that many other large families of unknotted satellite operators are injective) by studying the behavior of surface invariants under satellite operations. 

 In Section \ref{sec: disk maps and doubling}, we prove an analogue of \cite[Corollary 1.5]{hedden_whitehead_doubles} showing that if a pair of slice disks are distinguished by knot Floer homology then their positive Whitehead doubles will be distinguished as well. Given a slice disk $D$ of $K$, we denote the induced element $F_D(1)$ in $\widehat{HFK}(S^3,K)$ by $t_D$.

\begin{mainthm}\label{thm:1}
    %Suppose a knot $K \subset S^3$ bounds slice disks $D_1, D_2 \subset B^4$. If $t_{D_1} \ne t_{D_2}$, then $t_{\Wh(D_1)} \ne t_{\Wh(D_2)}$.
    %Suppose $K \subset S^3$ bounds slice disks $D_1, D_2 \subset B^4$. If $t_{D_1} \ne t_{D_2}$, then $t_{\Wh(D_1)} \ne t_{\Wh(D_2)}$.
%Let $D_1$ and $D_2$ be slice disks with $\partial D_1=\partial D_2$. If $t_{D_1} \ne t_{D_2}$, then $t_{\Wh(D_1)} \ne t_{\Wh(D_2)}$ as well. 
   Let $K$ be a knot in $S^3$ with slice disks $D_1$ and $D_2$. If $t_{D_1} \ne t_{D_2}$, then $t_{\Wh(D_1)} \ne t_{\Wh(D_2)}$ as well.
\end{mainthm}

In fact, we prove an analogous result for many satellite operations, including positive $(p,1)$-cables (Theorem~\ref{thm:1-1}), the Mazur pattern (Example~\ref{ex:mazur}), and a large class of satellite patterns introduced by Levine  \cite{levine_doubling_operators} that generalizes Whitehead doubling (Proposition~\ref{prop: doubled invariants}).  %We prove this statement for a much larger class of satellites introduced by Levine \cite{levine_doubling_operators}, which generalizes Whitehead doubling (as well as some other common unknotted patterns). 
Our proof relies on an analysis of the type A module associated to the pattern knot embedded in the solid torus. In particular, our proof holds for any unknotted pattern satisfying a technical condition on the type A structure, which we expect to hold for many unknotted patterns. % We have restricted our analysis to doubling and cabling operators because these patterns are extremely useful in the analysis and construction of exotic phenomena in the four-ball (see Section \ref{sec: exotically knotted surfaces}. But, we expect the results in this section to hold for many unknotted patterns. 

Our corresponding result for cabling patterns shows that $(p,1)$-cabling does more than preserve the difference between disks distinguished by knot Floer homology --- it increases the distance between them. To state this precisely, we recall the operation of \emph{(internal) stabilization} of a knotted surface as described in \cite[\S2.1]{baykur_sunukjian_knotted_surf_stab}: given an embedded surface $S \subset B^4$, choose an embedded 3-dimensional 1-handle $h \approx [-1,1]\times D^2$ in $B^4$ that intersects $S$ only along $\{\pm1\} \times D^2$, then increase the genus of $S$ by removing $\{\pm1\}\times D^2$ from $S$ and gluing in $[-1,1]\times D^2$. We show that the $(p,1)$-cables of disks which are distinguished by knot Floer homology must always have stabilization distance at least $p$. 

\begin{mainthm}\label{thm:1-1}
Let $K$ be a knot in $S^3$ with slice disks $D_1$ and $D_2$. If $t_{D_1} \ne t_{D_2}$, then for any integer $p\ge 1$, we have $t_{(D_1)_{p,1}} \ne t_{(D_2)_{p,1}}$. Furthermore, the stabilization distance between $(D_1)_{p,1}$ and $(D_2)_{p,1}$ is at least $p$.
\end{mainthm}

Below, we describe applications of these results and related techniques to the study of exotic surfaces. 

\subsection{Exotically knotted surfaces}\label{sec: exotically knotted surfaces}

%While the results from Section \ref{sec:injective satellite} %The  show that smoothly distinct disks often remain distinct under various satellite operations; things are quite different in the topological category. 
While the results above show that smoothly distinct disks often remain distinct under various satellite operations, things are quite different in the topological category.  Applying the work of Conway and Powell \cite[Theorem 1.2]{conway2021embedded}, we show that for any unknotted, winding number zero pattern $P$ and any pair of disks $D_1,D_2 \subset B^4$ with the same boundary, the satellite disks $P(D_1)$ and $P(D_2)$ are topologically isotopic rel boundary. In particular, given any pair of slice disks for a knot $K \sub S^3$, their Whitehead doubles are always topologically isotopic. Therefore, the obstructive results from Section \ref{sec:injective satellite} imply the existence of vast numbers of exotic disks; by starting with disks which can be distinguished by knot Floer homology we can produce an exotic pair by applying the Whitehead doubling operator.

\addtocounter{mainthm}{-1}
\begin{maincor}\label{cor:thm1}
    Let $K$ be a knot in $S^3$ with slice disks $D_1$ and $D_2$. If $t_{D_1} \ne t_{D_2}$, then $\Wh(D_1)$ and $\Wh(D_2)$ are exotically knotted.
\end{maincor}

For example, we show that any nontrivial knot of the form $K\smallsum K\smallsum -K\smallsum -K$ (which we will write as $2K\smallsum -2K$) has a pair of slice disks distinguished by knot Floer homology, implying the following:

%\begin{theorem}\label{cor:1}
\begin{maincor}\label{cor:1}
If $K$ is a nontrivial knot, then the knot $\Wh(2K\smallsum -2K)$ % $\Wh(K\smallsum K\smallsum -K\smallsum -K)$ 
bounds an exotic pair of  disks.
\end{maincor}
%\end{theorem}
Recall that disks are called \emph{$n$-stably exotic} if they remain exotic after applying $n$ internal stabilizations.   Combining Theorem \ref{thm:1-1} with Corollary \ref{cor:1} yields many $n$-stably exotic pairs for arbitrarily large values of $n$, significantly expanding on the results of \cite{guth_one_not_enough_exotic_surfaces}.

As another example, we consider the infinite family of slice disks for $4_1 \smallsum 4_1$ obtained by roll-spinning \cite{fox:rolling} (see also Figure~\ref{fig:swallow-follow}); these are distinct up to topological isotopy rel boundary (e.g., by comparing their \emph{peripheral maps} \cite[\S3]{juhász_zemke_dist_slice_disks}).  In \cite{juhasz_zemke_stabilization_bounds}, Juh\'asz and Zemke  show that knot Floer homology is able to distinguish the first two members of this family of  disks (up to smooth isotopy rel boundary), hence their Whitehead doubles are exotic by Corollary~\ref{cor:thm1}. However, in its current state, knot Floer homology cannot distinguish all members of this infinite family of  disks. This would require naturality and functoriality of the link cobordism maps for integral coefficients as well as a formula for the basepoint moving action with intergal coefficients \cite{sarkar2015moving,Zemke2016QuasistabilizationAB}. Nevertheless, we are able to distinguish the entire family of slice disks for $\Wh(4_1 \smallsum 4_1)$ indirectly using Seiberg-Witten invariants, building on arguments of Gompf \cite{gompf:infinite,gompf:infinite-handle} and Akbulut \cite{akbulut:infinite}.

%In \cite{juhasz_zemke_stabilization_bounds}, Juh\'asz and Zemke consider an infinite family of slice disks (obtained by roll-spinning) bounded by the knot $4_1 \# 4_1$. These disks are known to be topologically inequivalent rel boundary, but Juh\'asz and Zemke also show that... \KH{[have to pause for now]}...d bounds an infinite family of disks that are pairwise nonisotopic (rel boundary), pairwise distinct disks (distinct topologically as well as smoothly) which are obtained by roll-spinning. However, in its current state, knot Floer homology cannot establish the existence of the infinite family of pairwise exotic slice disks for $\Wh(4_1\#4_1)$ since the naturality and functoriality of the link cobordism maps have not been established for integral coefficients. Instead, we prove the result using Seiberg-Witten invariants.
%In \cite{juhász_zemke_dist_slice_disks}, Juh\'asz and Zemke  implicitly show that the knot $4_1 \# 4_1$ bounds an infinite family of disks that are pairwise nonisotopic (rel boundary), pairwise distinct disks (distinct topologically as well as smoothly) which are obtained by roll-spinning. However, in its current state, knot Floer homology cannot establish the existence of the infinite family of pairwise exotic slice disks for $\Wh(4_1\#4_1)$ since the naturality and functoriality of the link cobordism maps have not been established for integral coefficients. Instead, we prove the result using Seiberg-Witten invariants.

%Finally, we produce an infinite family of pairwise exotic slice disks in the four-ball. 
\addtocounter{mainthm}{1}
\begin{mainthm}\label{thm:2}
    The knot $\Wh(4_1 \# 4_1)$ bounds an infinite family of pairwise exotic slice disks.
\end{mainthm}

Applying techniques of Akbulut and Ruberman \cite{akbulut-ruberman,ruberman}, we can upgrade the examples in Theorem~\ref{thm:2} to drop the ``rel boundary'' condition and produce an infinite collection of exotic slice disks that are \emph{absolutely exotic}, i.e., not related by any smooth isotopy of $B^4$.

\begin{maincor}\label{cor:absolute}
There exists a knot in $S^3$ that bounds an  infinite family of exotic slice disks in $B^4$ that are not related by any smooth isotopy of $B^4$.
\end{maincor}

%\subsection{Cabling and stabilization}

%Our techniques extend to generate examples of stably exotic disks in the four-ball as well. Recall that disks are called \emph{$n$-stably exotic} if they remain exotic after applying $n$ internal stabilizations.  We show that the $(p,1)$-cables of disks which are distinguished by knot Floer homology must always have stabilization distance at least $p$. 

%\begin{mainthm}\label{thm:1-1}
%Let $K$ be a knot in $S^3$ with slice disks $D_1$ and $D_2$. If $t_{D_1} \ne t_{D_2}$, then for any integer $p\ge 1$, we have $t_{(D_1)_{p,1}} \ne t_{(D_2)_{p,1}}$. Furthermore, the stabilization distance between $(D_1)_{p,1}$ and $(D_2)_{p,1}$ is at least $p$.
%\end{mainthm}

%Of course, combining Theorem \ref{thm:1-1} with Corollary \ref{cor:1} yields many exotic pairs with arbitrarily large stabilization distance.

There is a close connection between exotically knotted surfaces and exotic 4-manifolds, especially through branched covering and surgery operations. For example, Finashin-Kreck-Viro constructed the first examples of exotically unknotted (nonorientable) closed surfaces in $S^4$ by proving that their branched double covers are exotic \cite{finashin1987exotic}. Furthermore, it has been hoped that an exotic $S^4$ might be realized through surgery, i.e. Gluck twist \cite{gluck1962embedding}, or taking a branched cover, along knotted 2-spheres in $S^4$.  By taking branched covers of   exotic surfaces that remain distinct after internal stabilization with $T^2$ (such as those in \cite{guth_one_not_enough_exotic_surfaces}), one obtains  natural candidates for exotic 4-manifolds that remain distinct after an external stabilization with $S^2 \! \times \! S^2$ (c.f., \cite{baykur_sunukjian_knotted_surf_stab}). However, we show that the branched covers of many stably exotic surfaces fail to yield stably exotic 4-manifolds:

\begin{mainthm}\label{thm:dbc}
There are surfaces $S,S' \subset B^4$ that are exotically knotted (rel boundary) yet whose branched double covers $\Sigma_2(B^4,S)$ and $\Sigma_2(B^4,S')$ are diffeomorphic (rel boundary). Moreover, we may construct $S,S' \subset B^4$ by internally stabilizing a pair of exotic slice disks $D,D' \subset B^4$ whose branched double covers $\Sigma_2(B^4,D)$ and $\Sigma_2(B^4,D')$ are exotic (rel boundary).
\end{mainthm}

Experts may notice that the surfaces in our proof of  Theorem~\ref{thm:dbc} are precisely those arising in \cite{guth_one_not_enough_exotic_surfaces}, which constituted the first known examples of stably exotic slice disks.

%Though we do not pursue it here, similar results can be obtained in the context of surgery operations. For example, let $S$ and $S'$ be the exotic genus-one surfaces from Theorem~\ref{thm:dbc}.  in (e.g., by considering log transforms along tori obtained by capping off once-stabilized

%Recently, the Khovanov invariants of surfaces have been successfully used to detect exotic surfaces \cite{Hayden_Sundberg_Kh_exotic_disks} and also to study Whitehead doubles of Seifert surfaces \cite{Hayden_Kim_Miller_Park_Sundberg_Kh_exotic_seifert}. Though Khovanov homology lacks the cut-and-paste techniques of Heegaard Floer homology, we observe that, in practice, the existence of a distinguishing pair of cycles for the companion disks can often be used to find a distinguishing pair of cycles for the satellite disks (see Section 7.)

\subsection{Khovanov homology} Several recent results have employed the cobordism maps on Khovanov homology to detect exotically knotted surfaces in $B^4$ \cite{Hayden_Sundberg_Kh_exotic_disks,Hayden_Kim_Miller_Park_Sundberg_Kh_exotic_seifert,hayden:atomic,lipshitz-sarkar:mixed}. In some contexts, the combinatorial nature of these invariants make their calculation more accessible than for their Floer-theoretic counterparts, and they appear to be better-suited to surfaces of higher genus (including those of non-minimal genus \cite{Hayden_Kim_Miller_Park_Sundberg_Kh_exotic_seifert}). Although Khovanov homology lacks the well-developed bordered theory of Heegaard Floer homology, some recent results have shed light on the effect of certain satellite operations on Rasmussen's invariant \cite{lewark-zibrowius} (see also \cite{zibrowius:2-tangles}) and the Khovanov cobordism maps \cite{Hayden_Kim_Miller_Park_Sundberg_Kh_exotic_seifert}, including distinguishing Whitehead doubles of certain Seifert surfaces. Here we extend the latter approach to the study of slice disks. Although we do not prove a general result analogous to Theorems~\ref{thm:1} and \ref{thm:1-1} for Khovanov homology, it can be shown that in many cases, the existence of cycles in Khovanov homology that distinguish a pair of slice disks can often be used to find cycles that distinguish satellites of those disks. We illustrate this with the following result.

\addtocounter{mainthm}{1}
\begin{mainprop}\label{prop:khovanov}
The knot $m(9_{46})$ bounds a pair of slice disks $D,D' \subset B^4$ whose Whitehead doubles $\Wh(D),\Wh(D')$ and $(2,1)$-cables $D_{2,1},D'_{2,1}$ are distinguished by their maps on Khovanov homology, $\Kh(\Wh(D))\neq \Kh(\Wh(D'))$ and $\Kh(D_{2,1})\neq\Kh(D'_{2,1})$. Moreover, the  maps $\operatorname{BN}(D_{2,1})$ and $\operatorname{BN}(D'_{2,1})$ on  Bar-Natan  homology over $\F_2[U]$ remain distinct after multiplication by $U$, hence these disks induce distinct maps  after one internal stabilization.
%    Let $D$ and $D'$ denote the standard slice disks bounded by the knot $9_{46}$, one of which is depicted in \Cref{fig:946-doubles}. The Khovanov cobordism maps $\Kh(\Wh(D))$ and  $\Kh(\Wh(D'))$ are distinct, as are the maps $\Kh(D_{2,1})$ and  $\Kh(D'_{2,1})$ induced by the $(2,1)$-cables of $D$ and $D'$. Moreover, over $\mathbb{F}_2[U]$, the Bar-Natan cobordism maps $\operatorname{BN}(D_{2,1})$ and $\operatorname{BN}(D'_{2,1})$ remain distinct even after multiplication by $U$, hence the once-stabilized disks induce distinct maps on Bar-Natan homology.
\end{mainprop}

It seems possible that further study of satellite operations on Khovanov homology and Bar-Natan homology may yield results analogous to Theorems~\ref{thm:1} and \ref{thm:1-1}.

%Experts may notice that a mild extension of the proof of Proposition~\ref{prop:khovanov} shows that...

\subsection{Conventions} Throughout this paper, the unknot will be denoted by $U$. The satellite of a knot $K$ with pattern $P$ will be written as $P(K)$. The notation $\nu(\, \cdot \,)$ will be used to denote an open tubular neighborhood, and $\bar \nu ( \, \cdot \, )$ its closure.  Isotopies of disks are always taken to be fixed on the boundary unless otherwise specified.

\section{Preliminaries}\label{sec: background}
In this section, we review the construction of satellite disks and establish a class of satellite patterns which are guaranteed to produce $\Z$-disks. 

\subsection{Satellite disks}\label{subsec: satellite disks}

Given a pattern $P$ with $P(U)=U$ and a slice disk $D$ which bounds a knot $K$, define the \emph{satellite disk} $P(D)$ as follows. Let $B$ be a neighborhood of a point in $D$. Then $D_0 :=D\smallsetminus B$ is a concordance in $B^4 \smallsetminus B \cong S^3 \times I$ between $U$ and $K$. We can then apply the pattern $P$, which gives a concordance $P(D_0)$ between $P(U)=U$ and $P(K)$. By capping off $P(U)=U$ with the trivial slice disk for the unknot, we obtain the satellite disk $P(D)$ for $K$.

A slice disk $D$ is called a \emph{$\Z$-disk} if the fundamental group of its complement is isomorphic to $\Z$. By the work of Conway and Powell, any two $\Z$-disks with common boundary are topologically isotopic rel.\ boundary \cite[Theorem~1.2]{conway_powell_characterisation_homotopy_ribbon}. We can arrange to work in this situation by choosing appropriate satellite patterns. Recall that the \emph{winding number} of a pattern $P$ is the algebraic intersection of $P$ with a generic meridional disk of the solid torus containing $P$. 

\begin{prop}\label{prop:whitehead doubles are Z-disks}
If $P$ is a winding number zero pattern with $P(U)=U$ and $D$ is a slice disk, then the satellite disk $P(D)$ is a $\Z$-disk.
\end{prop}
\begin{proof}
Choose a tubular neighborhood $\nu D$ of $D$. The satellite disk $P(D)$ is contained in $\nu D$, so we have a splitting
\[
B^4 \smallsetminus P(D) = (B^4 \smallsetminus \nu D) \cup (\bar \nu D\smallsetminus P(D)).
\]
Since homeomorphism class of $\bar\nu D \smallsetminus P(D)$ does not depend on the choice of $D$,
\[
\bar \nu D\smallsetminus P(D) \cong \bar \nu(D_0) \smallsetminus P(D_0),
\]
where $D_0$ denotes the trivial slice disk of an unknot $U$. It is clear that $P(D_0)$ is also a trivial slice disk of $P(U)=U$ and we may take $\bar \nu(D_0)$ to be the whole $B^4$. Therefore,  
\[
\bar \nu(D_0) \smallsetminus P(D_0) \cong B^4 \smallsetminus D_0 \cong D^2 \times (D^2 \smallsetminus \{pt\}),
\]
implying $\pi_1(\nu D \smallsetminus P(D))\cong \Z$.

The fundamental group of the intersection 
\[
(B^4 \smallsetminus \nu D)\cap (\bar \nu D\smallsetminus P(D)) \cong S^1 \times D,
\]
is clearly also $\Z$, and the natural maps
\[
\begin{split}
    \pi_1\left((B^4 \smallsetminus \nu D)\cap (\bar \nu D\smallsetminus P(D)\right) &\rightarrow \pi_1\left(B^4 \smallsetminus \nu D\right),\\
    \pi_1\left((B^4 \smallsetminus \nu D)\cap (\bar \nu D\smallsetminus P(D)\right) &\rightarrow \pi_1\left(\nu D \smallsetminus P(D)\right),
\end{split}
\]
are given by the inclusion of a meridional class of $\pi_1(B^4 \smallsetminus \nu D)$ and the multiplication map $\Z\xrightarrow{\times w(P)}\Z$, respectively, where $w(P)$ denotes the winding number of $P$. Since $P$ has winding number $0$ and $\pi_1(B^4 \smallsetminus \nu D)$ is normally generated by a meridian of $D$, we see
\[
\pi_1(B^4 \smallsetminus P(D)) \cong \pi_1(\bar \nu D \smallsetminus P(D)) \cong \Z.
\]
Therefore $P(D)$ is a $\Z$-disk.
\end{proof}

\begin{example}\label{ex:946-wh}
    Let $D$ denote the slice disk for $m(9_{46})$ from Figure~\ref{fig:946-doubles} and let $D'$ be obtained by applying the involution shown in the top left of Figure~\ref{fig:946-doubles}. By Proposition~\ref{prop:whitehead doubles are Z-disks}, the Whitehead doubled disks $\Wh(D)$ and $\Wh(D')$ are $\Z$-disks with the same boundary $\Wh(m(9_{46}))$, hence they are topologically isotopic rel.~boundary by \cite[Theorem~1.2]{conway_powell_characterisation_homotopy_ribbon}. In later sections we will prove that these disks are \emph{not} smoothly isotopic rel.~boundary.
\end{example}

\begin{figure}[hbt]
    \centering
    \includegraphics[scale=.4]{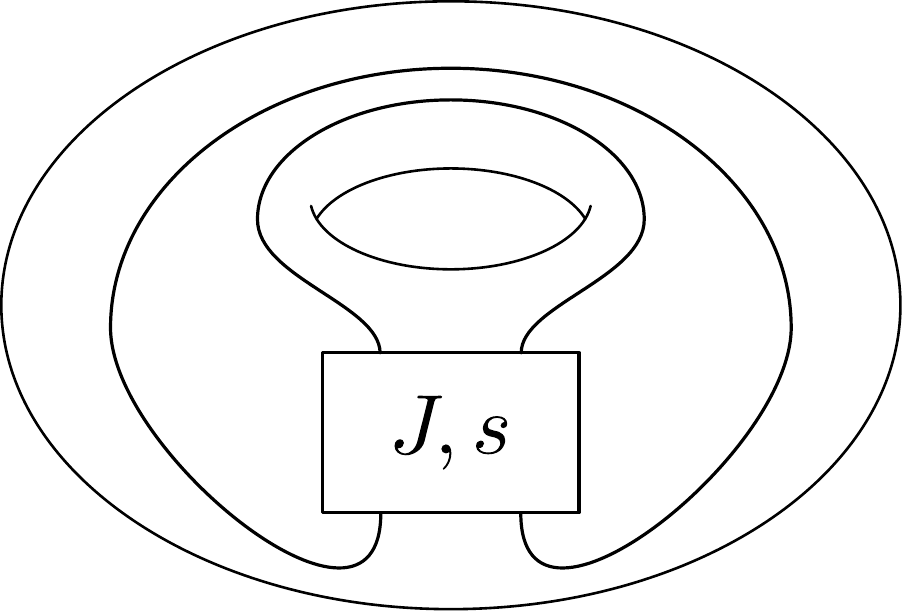}
    \caption{The doubling pattern associated to a knot $J$ with framing $s$.}
    \label{fig:doubling pattern}
\end{figure}

Knot doubling (which we now describe) produces an infinite family of patterns which satisfy the assumptions of \Cref{prop:whitehead doubles are Z-disks}. Given a knot $J$ and an integer $s$, the doubling pattern $\Lev_{J,s}$ is shown in \Cref{fig:doubling pattern}; the box labeled $J, s$ represents $s$-framed parallel copies of $J\smallsetminus \pt$. This operation generalizes Whitehead doubling, which is recovered by taking $J$ to be the unknot and $s = \pm1$. These satellites are described more carefully as twisted infections along components of the Borromean rings in \cite[Section 1]{levine_doubling_operators}. Since these patterns are unknotted, for any slice disk $D$, the doubled disk $\Lev_{J,s}(D)$ is a $\Z$-disk. Moreover, by the formula for the Alexander polynomial of a satellite knot, the knot $\Lev_{J, s}(K)$ is an Alexander polynomial 1 knot.

It is worth observing here that, in contrast with Proposition~\ref{prop:whitehead doubles are Z-disks}, patterns with nonzero winding number preserve most of the information about the fundamental group of a slice disk's complement.

\begin{proposition}\label{prop:nonzero-winding}
Let $P$ be a pattern with nonzero winding number.
\begin{enumerate}
    \item [\normalfont\textbf{(a)}] For any knot $K\subset S^3$, the inclusion-induced map $\pi_1(S^3 \setminus K) \cong \pi_1(S^3 \setminus \nu K) \to \pi_1(S^3 \setminus P(K))$ is injective.

        \item [\normalfont\textbf{(b)}] For any slice disk $D\subset B^4$, the inclusion-induced map $\pi_1(B^4 \setminus D) \cong \pi_1(B^4 \setminus \nu D)  \to \pi_1(B^4 \setminus P(D))$ is injective.
\end{enumerate}
%    Let $P$ be a pattern with $P(U)=U$ and nonzero winding number. For any slice disk $D \subset B^4$, the inclusion-induced map $\pi_1(B^4 \setminus \nu D) \to \pi_1(B^4 \setminus P(D))$ is injective.
\end{proposition}

\begin{proof}
    %\KH{An analysis similar to the one above shows that the inclusion-induced map $\pi_1(B^4\setminus \nu D) \to \pi_1(B^4 \setminus P(D))$ is injective, hence $\pi_1(B^4 \setminus D)$ naturally embeds as a subgroup of $\pi_1(B^4 \setminus P(D))$ if $w(P)\neq0$.}
We begin by recalling a key group-theoretic fact: Given groups $H$, $G$, and $G'$ and homomorphisms $f: H \to G$ and $f': H \to G'$, consider the amalgamated product $G *_H G'$. If $f$ and $f'$ are injective, then the natural maps of $G$ and $G'$ into $G *_H G'$ are injective (cf \cite[\S11]{rotman}).

For the claim in (a), we may assume that the pattern $P$ is not the identity pattern $S^1 \times \{pt\} \subset S^1 \times D^2$, as that case is trivial. The boundary $\partial \bar \nu K$ of the closed tubular neighborhood $\bar \nu K$ is an incompressible torus that separates $S^3 \setminus P(K)$ into $S^3 \setminus \nu K$ and $\bar \nu K \setminus P(K)$. This gives rise to a decomposition 
$$\pi_1(S^3 \setminus P(K)) \cong \pi_1(S^3 \setminus \nu K) *_{\footnotesize\raisebox{-2pt}{$\pi_1(\partial \bar \nu K)$}} \pi_1(\bar \nu K \setminus P(K)).$$
Since $\bar \nu K$ is incompressible, the maps on $\pi_1$ induced by including  $\partial \bar \nu K$ into $S^3 \setminus \nu K$ and $\bar \nu K \setminus P(K)$ are injective. It follows that the map $\pi_1(S^3 \setminus K) \cong \pi_1(S^3 \setminus \nu K) \to \pi_1(S^3 \setminus P(K))$ is also injective. 

The argument for (b) is even simpler. Consider the decomposition $$B^4 \setminus P(D) = \big(B^4 \setminus \nu D\big) \cup \big(\bar \nu D \setminus P(D)\big).$$ As in the proof of Proposition~\ref{prop:whitehead doubles are Z-disks}, the intersection of these two subsets has $\pi_1 \cong \Z$, generated by a meridian $\mu$ of $D$. The inclusion-induced map $\Z \to \pi_1(B^4 \setminus \nu D)$ is clearly injective, and the fact that $P$ has nonzero winding number implies that the other map $\Z \to \pi_1(\bar \nu D \setminus P(D))$ is also injective. As above, we conclude that the map from $\pi_1(B^4 \setminus D) \cong \pi_1(B^4 \setminus \nu D)$ into the amalgamated product $\pi_1(B^4 \setminus P(D))$ is injective.
\end{proof}

 Recall that the  \emph{peripheral map} of a slice disk $D \subset B^4$ is the inclusion-induced map $\pi_1(S^3 \setminus \partial D)\to \pi_1(B^4 \setminus D)$ (cf \cite[Definition 3.9]{juhász_zemke_dist_slice_disks}). Consider the following commutative diagram of peripheral maps.
\begin{center}
\begin{tikzcd}[row sep = small]
& \pi_1(B^4 \setminus \nu D) \ar[r] &  \pi_1(B^4 \setminus P(D))
\\
\pi_1(S^3 \setminus K) \cong \pi_1(S^3 \setminus \nu K)
 \ar[start anchor={[xshift=3.25ex]}, ur] \ar[start anchor={[xshift=4ex]}, dr] 
 \ar[r] & \pi_1(S^3 \setminus P(K))
 \ar[ur] \ar[dr] 		\\
& \pi_1(B^4 \setminus \nu D') \ar[r] & \pi_1(B^4 \setminus P(D'))
\end{tikzcd}
\end{center}
When $P$ has nonzero winding number, the horizontal arrows are injective, yielding the following corollary of Proposition~\ref{prop:nonzero-winding}.
 
\begin{corollary}\label{cor:nonzero-winding}
     Let $K \subset S^3$ be a knot that bounds slice disks $D,D' \subset B^4$ and let $P$ be a pattern with $P(U)=U$ and $w(P)\neq 0$. If the disks $D$ and $D'$ are distinguished by the kernels of their peripheral maps, then so are the satellite disks $P(D)$ and $P(D')$.
\end{corollary}

\begin{example}\label{ex:946-nonzero}
    Continuing the notation from Example~\ref{ex:946-wh}, we recall the well-known fact that the slice disks $D$ and $D'$ bounded by $m(9_{46})$ are distinguished by their peripheral maps' kernels. To see this, consider Figure~\ref{fig:946-exterior}, which exhibits the exterior of $D \subset B^4$. The diagram in part (a) arises from applying the procedure in \cite[\S6.2]{GS_4mflds} to obtain the exterior of $D\subset B^4$, where $\alpha$ is a decorative curve to be used later. Part (b) is obtained by isotopy, and part (c) by sliding one 1-handle over the other. After an isotopy that simplifies the diagram, we obtain (d); in this final diagram, the black and blue curves can be viewed as a handle diagram for $B^4$ in which $D$ is an unknotted disk bounded by the red curve.

\begin{figure}\center
\def\svgwidth{\linewidth}\input{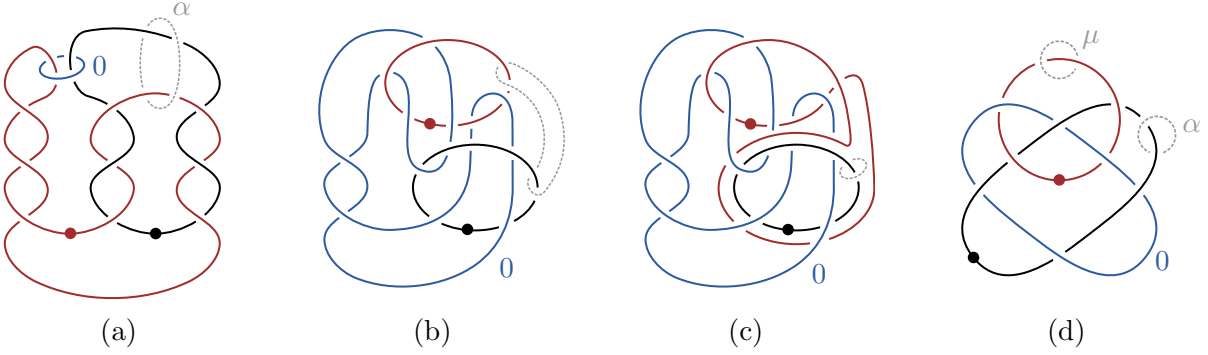}
\caption{Handle diagrams for the exterior of the slice disk $D\subset B^4$ bounded by $m(9_{46})$.}\label{fig:946-exterior}
\end{figure}

    To distinguish the disks $D$ and $D'$, observe that $\alpha$ is a curve in the knot complement $S^3 \setminus m(9_{46})$ that bounds a smoothly embedded disk in $B^4 \setminus D'$. On the other hand, we claim it is nontrivial in $\pi_1(B^4 \setminus D)$. To see this, choose based representatives of the loops $\mu$ and $\alpha$ in Figure~\ref{fig:946-exterior}(d). It is straightforward to check that $\pi_1(B^4\setminus D) \cong \langle \mu,\alpha : \alpha \mu \alpha^{-2} \mu^{-1} \rangle$. This group admits a surjection onto the symmetric group $S_3$ defined by sending $\mu$ to $(2\ 3)$ and $\alpha$ to $(1 \ 2 \ 3 )$ (cf \cite[Proposition 4.4]{akmr:stable}). It follows that $\alpha$ corresponds to a nontrivial element in $\pi_1(B^4\setminus D)$, yet it is trivial in $\pi_1(B^4 \setminus D')$ because it bounds a disk in $B^4 \setminus D'$. This implies that the peripheral maps associated to $D$ and $D'$ have distinct kernels, hence these disks are not  topologically isotopic rel boundary. By Corollary~\ref{cor:nonzero-winding}, it follows that $P(D)$ and $P(D')$ are not isotopic rel boundary for any unknotted pattern $P$ with nonzero winding number (including, for example, any $(p,1)$-cable with $p>1$).
\end{example}

\section{Floer Theoretic Background}

\subsection{Knot Floer homology and cobordism maps} \label{sec: knot Floer}

Knot Floer homology, defined in \cite{os_knotinvts} and independently in \cite{rasmussen_knotcompl}, is an invariant of knots, links, and cobordisms between them. In the terminology of \cite{zemke_linkcob}, link Floer homology associates to a decorated link $\mathbf{L}$ a curved complex  $CFL(S^3,\mathbf{L})$ over the commutative ring freely generated over $\mathbb{F}_2$ by formal variables which correspond to basepoints on the decorated link; to a decorated cobordism it associates a chain map between the curved complexes assigned to its ends. With the aim of being self-contained, we give the definitions of decorated links and cobordisms as well as the category of curved complexes here. For simplicity, we will write ``decorated links'' rather than \emph{multi-based links with a coloring} as in \cite{zemke_linkcob}.

\begin{defn}
A \emph{decorated link} is a 5-tuple $\mathbf{L}=(L,\mathbf{w},\mathbf{z},P,\sigma)$ which satisfies the following conditions:
\begin{itemize}
    \item $L$ is a link in $S^3$
    \item $\mathbf{w}$ and $\mathbf{z}$ are finite collections of basepoints in $L$, which intersect nontrivially with every component of $L$, such that closures of every component of $L\setminus (\mathbf{w}\cup \mathbf{z})$ has one endpoint in $\mathbf{w}$ and the other endpoint in $\mathbf{z}$
    \item $P$ is a finite set, and $\sigma$ is a function from $\mathbf{z}\cup\mathbf{w}$ to $P$, which defines a ``coloring'' of the basepoints.
\end{itemize}
\end{defn}

% , for $\mathbf{L}_i = (L_i,\mathbf{z}_i,\mathbf{w}_i,P,\sigma_i)$ with $i=0,1$

\begin{defn}
For $i=0,1$, let $\mathbf{L}_i = (L_i,\mathbf{z}_i,\mathbf{w}_i,P,\sigma_i)$ be decorated links. A \emph{decorated cobordism} between $\mathbf{L}_1$ and $\mathbf{L_2}$ is a 5-tuple $\mathbf{S}=(S,\Sigma_z,\Sigma_w,P,\sigma)$ which satisfies the following conditions.
\begin{itemize}
    \item $S$ is a properly embedded smooth surface in $S^3 \times [0,1]$, such that $\Sigma\cap (S^3 \times \{i\})=L_i$ for $i=0,1$.
    \item $\Sigma_w$ and $\Sigma_z$ are closed subsurfaces of $S$ such that $\Sigma_w \cup \Sigma_z = S$ and $\Sigma_w \cap \Sigma_z$ is a properly embedded 1-dimensional submanifold of $S$.
    \item $L_i \cap \Sigma_w = \mathbf{w}_i$, $L_i \cap \Sigma_z = \mathbf{z}_i$, and $L_i \cap \Sigma_w \cap \mathbf{w}_i$ and $L_i \cap \Sigma_z \cap \mathbf{z}_i$ consist of one point for $i=0,1$.
    \item $\sigma$ is a function from $\pi_0(\Sigma_w) \cup \pi_0(\Sigma_z)$ to $P$, such that every basepoint $p\in \mathbf{w}_i \cup \mathbf{z}_i$ which is contained in a connected component $S\in \pi_0(\Sigma_w)\cup \pi_0(\Sigma_z)$ satisfies $\sigma_i(p)=\sigma(S)$ for $i=0,1$.
\end{itemize}
\end{defn}

Having defined decorated links and cobordisms, we can describe the algebraic structure of knot Floer homology a bit more precisely. Given a finite set $P$, consider the commutative ring $R_P$, freely generated over the field $\mathbb{F}_2$ by the elements of $P$. The curved complex $CFL(S^3,\mathbf{L})$ is a bigraded module over $R_P$. A decorated link cobordism $\mathbf{S}$ between decorated links $\mathbf{L}_1$ and $\mathbf{L}_2$, induces an $R_P$-linear chain map $F_\mathbf{S}$ from $CFL(S^3,\mathbf{L}_1)$ to $CFL(S^3,\mathbf{L}_2)$, which we call the \emph{cobordism map} induced by $\mathbf{S}$.

In this paper, we will focus on knots, i.e. links with only one component (in which case we denote $CFL$ by $CFK$ instead), and concordances between them, i.e. cobordisms between knots which are homeomorphic to cylinders. We will endow all knots $K$ with the simplest possible decoration, consisting of two basepoints $\mathbf{w}=\{w\}$ and $\mathbf{z}=\{z\}$, two colors $P=\{w,z\}$ and the identity coloring function $\mathrm{id}:\{w\}\cup\{z\}\rightarrow \{w,z\}$. In this case, we have $R_P = \mathbb{F}_2[\U,\V]$ where $\U$ and $\V$ are the formal variables associated to the colors $w$ and $z$, and the curvature of the complex $CFK(S^3,\mathbf{K})$ vanishes, so that it becomes a chain complex over $\mathbb{F}_2[\U,\V]$. 

To make expressions simpler, we will abuse notation, and denote both the knot $K$ and the decorated knot $\mathbf{K} = (K, w, z)$ by $K$. We will denote by $CFK^\infty(S^3,K)$ the chain complex over $\mathbb{F}_2[\U,\U^{-1}]$ obtained from $CFK(S^3,K)$ by truncating with $\V=0$ and then localizing by adjoining the formal inverse of $\U$. Furthermore, we denote the truncation of $CFK(S^3,K)$ with $\U=\V=0$ by $\widehat{CFK}(S^3,K)$, and its homology by $\widehat{HFK}(S^3,K)$

In a similar way, a concordance between knots can be endowed with the simplest possible decoration, so $\Sigma_w$ and $\Sigma_z$ are both rectangles. The cobordism map
\[
F_\mathbf{S} : CFK(S^3,K_1)\rightarrow CFK(S^3,K_2),
\]
induced by the decorated cobordism $\mathbf{S}$, is $\mathbb{F}_2[\U,\V]$-linear, bidegree-preserving, and induces a quasi-isomorphism on $CFK^\infty$ by \cite[Theorem C]{zemke_linkcob}. Again, for simplicity, we will simply write $S$ for the decorated concordance described above. 

In general, the chain map $F_S$ is not well-defined if a decoration is not specified; due to the ambiguity induced by full rotations along knots, its homotopy class is defined uniquely only up to composition by the basepoint moving map, which was studied by Sarkar in \cite{sarkar2015moving}, of $K_1$ (or equivalently, $K_2$). However, when either $K_1$ or $K_2$ is unknotted (which is the context within which we will be working) the map does not depend on the choice of decoration, since the basepoint moving map on the unknot is the identity.

In particular, if $S$ is the concordance obtained by puncturing a slice disk $D$ of a knot $K$, the chain map
\[
F_S:CFK(S^3, U)\rightarrow CFK(S^3,K)
\]
is well-defined up to homotopy. Truncating by $\U=\V=0$ and then taking homology gives a cobordism map
\[
\widehat{F}_S:\widehat{HFK}(S^3, U)\rightarrow \widehat{HFK}(S^3,K),
\]
which is an invariant of the surface. Since $\widehat{HFK}(S^3, U)\simeq \mathbb{F}_2$, this map can be represented by the homology class
\[
t_D \coloneqq \hat{F}_S(1) \in \widehat{HFK}(S^3,K).
\]
Note that the class $t_D$ depends only on the smooth isotopy class of the given slice disk $D$. This fact was first observed in \cite{juhasz2016concordance}.

\subsection{Bordered Floer homology}\label{sec: bordered backgroun}

Bordered Floer homology is a package of invariants of 3-manifolds with boundary \cite{LOT_bordered_HF}. It associates to a surface $F$ a differential graded algebra $\cA(F)$ and to a 3-manifold $Y$ parametrized boundary both a type D structure over $\cA(-F)$, called $\CFDh(Y)$, and type A structure over $\cA(F)$, called $\CFAh(Y)$. 

Recall that a type $D$ structure over $\cA$ is a pair $(N, \delta^1)$ where $N$ is a graded $\F$-module and
\[
\delta^1: N \ra (A \otimes N)[1],
\]
is an $\F$-module map satisfying 
\[
(\mu_2 \otimes \bI_N)\circ(\bI_A \otimes \delta^1) \circ \delta^1 + (\mu_1\otimes \bI_N) \circ \delta^1 = 0.
\]
We will frequently be interested in \emph{type D structure homomorphisms}, which are $\F$-module maps $f^1: N_1 \ra A \otimes N_2$ satisfying 
\[
(\mu_2 \otimes \bI_{N_2}) \circ (\bI_A \otimes f^1) \circ \delta_{N_1}^2 + (\mu_2 \otimes \bI_{N_2}) \circ (\bI_A \otimes \delta^1_{N_2})\circ f^1 + (\mu_1\otimes \bI_{N_2}) \circ f^1 = 0.
\]
The map $\delta^1$ can be iterated to define maps 
\[
\delta^k: N \ra (A^{\otimes k} \otimes N)[k],
\]
where $\delta^0 = \bI_{N}$ and $\delta^k = (\bI_{A^{\otimes(k-1)}} \otimes \delta^1) \circ \delta^{k-1}.$ Similarly, a type D homomorphism $f^1$ can be used to define maps 
\[
f^k: N_1 \ra (A^{\otimes k} \otimes N_2)[k-1], \;\; f^k(x) = \sum_{i+j=k-1}  (\bI_{A^{\otimes(i-1)}} \otimes \delta^j_{N_2}) \circ  (\bI_{A^{\otimes i}} \otimes f^1) \circ \delta^i_{N_2}.
\]

A type A structure is a graded $\F$-module $M$ equipped with operations 
\[
m_i: M \otimes A^{\otimes (i-1)} \ra M[2-i],
\]
satisfying 
\begin{align*}
    &\sum_{i+j = n+1} m_i(m_j(x \otimes a_1 \otimes \hdots \otimes a_{j-1}) \otimes \hdots \otimes a_{n-1}) \\
+ &\sum_{i+j=n+1}\sum_{\ell = 1}^{n-j} \mu_i(x, a_1 \otimes \hdots \otimes a_{\ell-1} \otimes \mu_j(a_\ell \otimes \hdots \otimes a_{\ell + j - 1}) \otimes\hdots \otimes a_{n-1}) = 0.
\end{align*}
Given both a type A and a type D structure, we define the box tensor product of $M$ and $N$ to be the $\F$-module $M \boxtimes N = M \otimes_\F N$, equipped with differential 
\[
\partial^\boxtimes(x \otimes y) = \sum_{k=0}^\infty (m_{k+1} \otimes \bI_N)(x \otimes \delta^k(y)).
\]
The differential on the box tensor product is often graphically depicted as 
\[
\partial^\boxtimes = 
\begin{tikzcd}
\, \ar[dd, dashed] & \, \ar[d, dashed] \\
& \delta \ar[dd, dashed] \ar[dl, Rightarrow] \\
m \ar[d, dashed] & \\
\, & \,
\end{tikzcd}.
\]
Given a map of type D structure $f^1: N_1 \ra N_2$, we will be interested in maps of the form $\bI_M \boxtimes f^1: M \boxtimes N_1 \ra M \boxtimes N_2$, which are defined as
\[
(\bI_M \boxtimes f^1)(x\otimes y) = \sum_{k=1}^\infty (m_{k+1} \otimes \bI_{N_2}) (x \otimes f^k(y)).
\]
Graphically, these map are given by
\[
\bI_M \boxtimes f^1 = 
\begin{tikzcd}
\, \ar[dddd, dashed] & \, \ar[d, dashed] \\
& \delta_{N_1} \ar[d, dashed] \ar[dddl, Rightarrow] \\
& f^1 \ar[d, dashed] \ar[ddl] \\
& \delta_{N_2} \ar[dd, dashed] \ar[dl, Rightarrow] \\
m \ar[d, dashed] & \\
\, & \,
\end{tikzcd}.
\]

The utility of these two algebraic objects is manifest in the bordered Floer homology pairing theorem \cite[Theorem 1.3]{LOT_bordered_HF}: if $Y_1$ and $Y_2$ are three manifolds with $\partial Y_1 \cong F \cong \partial Y_2$, the hat-version of the Heegaard Floer homology of $Y_1 \cup_F Y_2$ is recovered by taking the box tensor product of the bordered type A and D structures associated to $Y_1$ and $Y_2$:
\[
	\CFh(Y_1\cup Y_2) \simeq \CFAh(Y_1) \boxtimes_{\cA(F)} \CFDh(Y_2).
\]
Bordered Floer theory also recovers (the $U = 0$ and hat versions of) knot Floer homology \cite[Theorem 11.21]{LOT_bordered_HF}: given a doubly pointed bordered Heegaard diagram $(\mathcal{H}_1, w, z)$ for $(Y_1, \partial F, K)$ and a bordered Heegaard diagram $(\mathcal{H}_2, z)$ with $\partial Y_1 \cong F \cong -\partial Y_2$, then 
\[
	\HFKh(Y_1 \cup Y_2, K) \cong H_*(\CFAh(\mathcal{H}_1, w, z) \boxtimes_{\cA(F)} \CFDh(\mathcal{H}_2, z)).
\]
The analogous version of the pairing theorem for $\HFK^-$ will be relevant to us as well. The pairing theorems give an effective way to study satellites: if $P(K)$ is the satellite of $K$ with pattern $P$, then $\HFKh(S^3, P(K))$ can be computed by finding a doubly pointed bordered Heegaard diagram $\mathcal{H}_P$ for the pattern $P$ in the solid torus and computing the box tensor product $\CFAh(\mathcal{H}_P) \boxtimes \CFDh(S^3 \smallsetminus K)$. 

Much in the same way, by passing through bordered Floer homology, we can compute link cobordism maps associated to satellite concordance in terms of the  maps associated to the original concordance.

\begin{thm}[{\cite[Theorem 2]{guth_one_not_enough_exotic_surfaces}}]\label{Theorem: satellite formula for concordance maps}
    Let $C: K \ra K'$ be a smooth concordance. Then, there exists a map $F: \CFDh(S^3 \smallsetminus K) \ra \CFDh(S^3 \smallsetminus K')$ induced by $C$, such that for any pattern knot $P$ in the solid torus and a suitable decoration on $P(C)$, the following diagram commutes up to homotopy:

    \begin{center}
        \begin{tikzcd}
            \CFAh(\mathcal{H}_P) \boxtimes \CFDh(S^3 \smallsetminus K) \ar[d, "\bI_{\CFAh(\mathcal{H}_P)} \boxtimes F"]  \ar[r, "\simeq"]  
                & \CFKh(S^3, P(K)) \ar[d, "F_{P(C)}"]  \\
            \CFAh(\mathcal{H}_P) \boxtimes \CFDh(S^3 \smallsetminus K') \ar[r, "\simeq"]
                & \CFKh(S^3, P(K^\prime)),
        \end{tikzcd}
    \end{center}
where $\mathcal{H}_P$ is a doubly pointed, bordered Heegaard diagram for $P\sub S^1\times D^2$, $P(K)$ and $P(K^{\prime})$ are satellites of $K$ and $K^{\prime}$, and $P(C)$ is the satellite concordance induced by $P$. The horizontal arrows are given by the pairing theorem \cite{LOT_bordered_HF}. The analogous statement for $\CFK^-$ holds as well.
\end{thm}

\section{Disk maps and satellite operators}\label{sec: disk maps and doubling}

In this section, we prove the key proposition which will lead to the proof of \Cref{thm:1}. 

\begin{prop}\label{prop: doubled invariants}
Let $K$ be a knot in $S^3$ and let $D_1$ and $D_2$ be a pair of slice disks for $K$. If the invariants $t_{D_1}$ and $t_{D_2}$ are distinct in $\HFKh(S^3, K)$, then for any knot $J$ and $s < 2\tau(J)$, $t_{\Lev_{J, s}(D_1)}$ and $t_{\Lev_{J,s}(D_2)}$ are distinct as well. 

\end{prop} 

We begin by giving a simple criterion for determining when an unknotted pattern $P$ is guaranteed to produce $\HFKh$-distinguishable disks in terms of the structure of $\CFAh(S^1\times D^2, P)$. 

When $P$ is an unknotted pattern, we have that
\[
\CFAh(S^1\times D^2, P) \boxtimes \CFDh(S^3\smallsetminus U) \simeq \CFKh(S^3, U).
\]
Therefore, there is some element $a$ in $\CFAh(S^1\times D^2, P)$ such that $a \otimes v$ generates homology of $\CFKh(S^3, U)$, where $v$ is the unique element in $\CFDh(S^3\setminus U)$. On homology, $\bI_P \boxtimes F$ is therefore determined by the image of $a\otimes v$. 

\begin{lemma}[No cancellation lemma]\label{lem: no cancellation} 
Let
\[
F: \CFDh(S^3\setminus U) \ra \CFDh(S^3 \setminus K),
\]
be a morphism of type D structures which tensors with the identity map of $\CFAh(S^1\times D^2, \lambda)$ to give a nontrivial map, where $\lambda$ is the knot $S^1 \times \pt \sub S^1 \times D^2$. Let $a$ be an element of $\CFAh(S^1 \times D^2, P)$ such that $a \otimes v$ generates the homology of $\CFAh(S^1 \times D^2, P)\boxtimes \CFDh(S^3\setminus U)$ and extend $\{a\}$ to a basis for $\CFAh(S^1 \times D^2, P)$. If the coefficient of $a$ is zero in every $\cA_\infty$ operation $m_k(b, \rho_{i_1}, \hdots, \rho_{i_{k-1}}) $ of $\CFAh(S^1 \times D^2, P)$ which preserves the filtration, then the map 
\[\bI_P \boxtimes F: H_*(\CFAh(S^1 \times D^2, P)\boxtimes \CFDh(S^3\setminus U)) \ra H_*(\CFAh(S^1 \times D^2, P)\boxtimes \CFDh(S^3\setminus K))\] 
is nontrivial.
\end{lemma}
\begin{proof}
Let $x = F(v)$, where $F$ is the map of type D structures induced by $F_C$. We can write 
\[
x = 1\cdot \theta + \sum_{I \in \{1, 2, 3, 12, 23, 123\}} \rho_I \theta_I,
\]
for some $\theta_I \in\CFDh(S^3\setminus K)$. The term $\theta$ must be nontrivial, since we have assumed that  $\bI_{\CFAh(S^1\times D^2, \lambda)} \boxtimes F$ has nontrivial image (this follows from the fact that $\CFAh(S^1\times D^2, \lambda)$ has no nontrivial $\cA_\infty$-operations and so all other terms in $F(x)$ are annihilated after taking the box tensor product.) $(\bI_P \boxtimes F)(a\otimes v)$ is defined to be
\begin{align*}
(\bI_P \boxtimes F)(a\otimes v) &= \sum_{k=1}^\infty (m_{k+1}\otimes \bI_P) \circ (a \otimes F^k(v))\\
&= a \otimes \theta + \text{ other terms}.
\end{align*}
The term $a\otimes \theta$ could be canceled if, for some $k$, $F^k(v)  = \rho_{i_1}\otimes \hdots\otimes \rho_{i_k} \otimes \theta$ and there is an operation of the form $m_{k+1}(a, \rho_{i_1}, \hdots, \rho_{i_k})$ in which $a$ appears with nonzero coefficient. However, we have assumed that no such operations exist. 

Moreover, when we pass to homology, there are no relations between $a\otimes \theta$ and any other elements of $H_*(\CFAh(S^1 \times D^2, P)\boxtimes \CFDh(S^3\setminus K))$, since $a\otimes \theta$ could only appear as a term in the boundary of another element if there were a filtration preserving operation of the form $m_{k+1}(b, \rho_{i_1}, \hdots, \rho_{i_k})$ in which $a$ appeared with nonzero coefficient. Again, no such operations exist. 

Therefore, $a \otimes \theta$ appears as a non-canceling term in the expansion of $(\bI_P \boxtimes F)(a\otimes v) \in H_*(\CFAh(S^1 \times D^2, P)\boxtimes \CFDh(S^3\smallsetminus K))$, from which it follows that $\bI_P \boxtimes F$ has nontrivial image.
\end{proof}

Therefore, to show knot Floer homology will continue to distinguish slice disks after applying an unknotted pattern $P$, it suffices to show the type A structure of the pattern satisfies the hypotheses of Lemma \ref{lem: no cancellation}. To demonstrate this process, we consider several concrete examples.

\begin{example}
    The type A structure for the positive Whitehead double is computed by Levine \cite{levine_doubling_operators} and is shown below: 

\begin{center}
\begin{tikzcd}[row sep = large, column sep= large]
c &
	c' \ar[l]&
		\\
&	
	b\ar[lu,"\rho_3\rho_2\rho_1"]\ar[d,"\rho_1"] &
		b' \ar[lu,"\rho_3\rho_2\rho_1"] \ar[l] \ar[ld,"\rho_{123}"] \ar[d,"\rho_1"] \ar[dd, bend left, "\rho_{12}"] \\
&
	a &
		a'\ar[l,"1+\rho_{23}"] \ar[d,"\rho_2"] \\
&
	&
		d \ar[ul,"\rho_3"].
\end{tikzcd}
\end{center}
In the diagram above, an arrow of the form $x \xra{\rho_{i_1} \hdots \rho_{i_k}} y$ indicates that $m_{k+1}(x, \rho_{i_1}, \hdots, \rho_{i_k}) = y$. Arrows pointing left lower the filtration. A short computation shows that $H_*(\CFAh(S^1\times D^2, \Wh) \boxtimes \CFDh(S^3\smallsetminus U)) = \F\langle b \otimes v \rangle$. By Lemma \ref{lem: no cancellation}, we can verify the nontriviality of $\bI_{\Wh} \boxtimes (F_1 + F_2)$ by simply analyzing those $\cA_\infty$ operations which output $b$, i.e. those arrows which point into $b$. Clearly, there is one such arrow, as $b$ appears in $m_1(b')$, but this operation lowers the filtration level. Therefore, by Lemma \ref{lem: no cancellation}, the map $\bI_{\Wh} \boxtimes (F_1 + F_2)$ must be nontrivial. 

\end{example}

\begin{example}\label{ex:mazur}
As another example, consider the Mazur pattern $M$. The type $A$ structure for the Mazur pattern was computed in \cite{petkova2021twisted}, and the $U = 0$ truncation is is shown below. 

\begin{center}
\begin{tikzcd}[row sep = large, column sep= large]
x_0 & 
x_1 \ar[l,"\rho_2"] & 
x_2 \ar[l,"\rho_1"]\ar[ll, bend right,"\rho_{12}" above]& 
x_3 \ar[ld,"\rho_2"]& 
x_4 \ar[l,"\rho_1"] \ar[lld, bend right=45,"\rho_{12}" above ]& x_5 & x_6 \\
& 
y_1 & 
y_2 & 
y_3 \ar[ll, bend left, "\rho_2\rho_1" above] & 
y_4 \ar[l,"\rho_1"]\ar[lll, bend left, "\rho_{12}\rho_1"]& 
y_5 & y_6
\end{tikzcd}
\end{center}
Here, $\iota_0\cdot\CFAh(S^1\times D^2, M) = \langle x_0, x_2, x_4, y_0, y_2, y_4 \rangle$ and the homology of $\CFAh(S^1\times D^2, M) \boxtimes \CFDh(S^3 \smallsetminus U)$ is generated by $y_4 \otimes v$. In this example, there are no arrows into $y_4$ at all. Therefore, $\HFKh$-distinguishable disks will remain distinguishable after applying the Mazur pattern.  
\end{example}

Of course, not all unknotted satellite patterns satisfy the hypotheses of the ``no cancellation lemma'', and this failure can give rise to disks which are distinguishable by knot Floer homology, but have satellites which are not. 

\begin{example}

    Consider our running example, $K = m(9_{46})$. The knot Floer homology calculator \cite{szabo_knot_floer_calc} computes that $\CFK(K)$ consists of a singleton $x$ and two unit boxes, generated by elements $a_i, b_i, c_i, e_i$ for $i \in \{1, 2\}$. 
    
    \[
    \begin{tikzcd}[column sep = large, row sep = large] b_i\ar[d,"V"] & a_i \ar[l, "U"] \ar[d, "V"]\\
    e_i & c_i \ar[l, "U"] 
    \end{tikzcd} \quad x \;\;
    \]

    The elements $x, e_i$ are in bigrading $(0,0)$. Since the two slice disks, $D$ and $D'$, for $K$ are symmetric, the work of Dai-Hedden-Mallick \cite{DaiHeddenMallick_corks_invol_floer} can be applied to show that 
    \[
    t_{D} + t_{D'} = e_1 + e_2 \in \CFKh(K).
    \]
    Hence, the two symmetric slice disks for $m(9_{46})$ can be distinguished by knot Floer homology. See \cite[Proposition 4.1]{guth_one_not_enough_exotic_surfaces} for a more detailed computation. \\
    
    This basis for $\CFK(K)$ induces a basis for $\CFDh(S^3\smallsetminus K)$ by \cite[Theorem 11.26]{LOT_bordered_HF}, which is shown below:

    \[
\begin{tikzcd}[row sep=1.1cm, column sep=1.1cm] b_i \ar[d, "\rho_1"]
            & y^1_i \ar[l, "\rho_2"]
                & a_i \ar[l, "\rho_3"] \ar[d, "\rho_1"]& \\
        y^2_i 
            &     
                & y^4_i 
                    & x  \ar[loop right, "\rho_{12}"] \\
        e_i \ar[u, "\rho_{123}"]
            & y^3_i \ar[l, "\rho_2"]
                & c_i \ar[l, "\rho_3"] \ar[u, "\rho_{123}"]
                    &
\end{tikzcd}
\]

By Theorem \ref{Theorem: satellite formula for concordance maps}, there is an induced map of type D structures (in fact, there is a unique map in this case), and it is given by 
\[
F:\CFDh(S^3\smallsetminus U) \ra \CFDh(S^3\smallsetminus K), \hspace{1cm} v \mapsto e + \rho_3 y^2 + \rho_1 y^3,
\]
where $w = w_1 + w_2$, for $w \in \{e, y^2, y^3\}.$ \\

Consider the $(2,-1)$ cabling pattern, which we denote $C_{2,-1}$. The type A structure is computed in \cite[Section 8]{OzStSz_concordanceHoms} and is shown below:
\[
\begin{tikzcd}
    & A_2 \ar[dl,"\rho_2"]\ar[dd,"\rho_{23}"] & A_1 \ar[l, "\rho_2\rho_1"] \ar[lld, bend right=45, "\rho_2\rho_{12}" above]\ar[ddl,"\rho_2\rho_{123}"]\\
    X\ar[dr,"\rho_3"] & & \\
    & B_2 & B_1
\end{tikzcd}
\]
The element $X$ is the sole basis element in $\iota_0 \cdot \CFAh(S^1 \times D^2, C_{2,-1})$. Therefore, we see that 
\[
(\bI_{\CFAh(S^1 \times D^2, C_{2,-1})} \boxtimes F)(X \otimes v) = X \otimes e + B_2 \otimes y^2.
\]
However, this element vanishes on homology. Note that 
\[
\partial^\boxtimes(A^2 \otimes y^3) = (m_2 \otimes \bI_{\CFDh(S^3\smallsetminus K)})(A^2 \otimes \delta^1(y^3)) = m_2(A^2 \otimes \rho_2)\otimes e = X \otimes e
\]
and 
\[
\partial^\boxtimes(A^1 \otimes y^3) = 
(m_3 \otimes\bI_{\CFDh(S^3\smallsetminus K)})(A^1, \delta^2(y^3))= m_3(A^1, \rho_2, \rho_{123}) \otimes y^2 = B_2 \otimes y^2.
\]
In fact, by considering the type A structure for $C_{p,-1}$ (see \cite[Section 8]{OzStSz_concordanceHoms}), we see that the same argument holds for any positive $p$. This example is in strong contrast to the behavior of positive cabling patterns; see Section \ref{sec: cables}.

\end{example}

Let us now turn to Levine's infinite family of doubling patterns. 

\begin{proof}[Proof of \Cref{prop: doubled invariants}]
The type A structure associated to the pattern $\Lev_{J, s}$ in the solid torus, which we denote $\CFAh(S^1 \times D^2,\Lev_{J, s})$, is computed by Levine. As the knot $J$ is unspecified, the computation is quite complicated. We will only need to consider a small subset of the $\cA_\infty$-operations, so we will only provide a terse description of this module, and refer the careful reader to \cite[Section 3]{levine_doubling_operators}. 

When $s < 2\tau(J)$, the type A structure $\CFAh(S^1 \times D^2,\Lev_{J, s})$ is generated by elements 
\[
    A^j, A'^j, B^j, B'^j, C^j, C'^j, D^j, D'^j, E^j_i, E'^j_i, F^j_i, F'^j_i, G^j_i, G'^j_i, H^j_i, H'^j_i, 
\]
where $0 \le j \le \text{rank}(\iota_0 \cdot \CFDh(S^3 \smallsetminus J))$ and $1 \le i \le k_j$, where $k_j$ is the length of the arrow between the $(2j-1)^{th}$ and $(2j)^{th}$ generators of $\iota_0 \cdot \CFDh(S^3 \smallsetminus J)$.

In this case (where $s < 2\tau(J)$), Levine shows that the homology of  $\CFAh(S^1 \times D^2,\Lev_{J, s}) \boxtimes \CFDh(S^3 \smallsetminus U)$ is generated by $A^0 \otimes v$, where $v$ is the single generator of $\CFDh(S^3 \smallsetminus U).$ Since we are only interested in the image of this element under the maps induced by the doubles of the disks in consideration, it suffices to consider those $\cA_\infty$-operations with outputs of the form \[A^0 + \text{other terms}.\] There is only one such operation, \[m_1(A'^0) =  A^0,\] but this operation does not preserve the filtration (see \cite[Figure 22]{levine_doubling_operators}. Basis elements are arranged in columns according to the Alexander filtration). Therefore, all criteria of \Cref{lem: no cancellation} are satisfied, so \Cref{prop: doubled invariants} follows.
\end{proof}

\subsection{Disks distinguished by knot Floer homology}

From \Cref{thm:1}, we see that to construct exotic disks bounding $\Lev_{J,s}(K)$, it suffices to construct two slice disks of $K$ which are distinguished by their induced maps on knot Floer homology. It is very easy to find such examples; a systematic way of finding them comes from a 4-dimensional construction called \emph{deformation-spinning}, which we now recall.

\begin{defn}\label{def:spinning}
Let $a$ be a properly embedded smooth arc in $D^3$.
  Furthermore, let $\phi \colon I \times D^3 \to D^3$ be an isotopy of $D^3$ such that
  $\phi_0 = \mathrm{id}_{D^3}$, $\phi_t \vert_{\partial D^3} = \mathrm{id}_{\partial D^3}$ for every $t \in I$, and $\phi_1(a) = a$.
  Then the \emph{deform-spun} slice disk $D_{a, \phi} \subset D^4$ is defined by taking
  \[
  \bigcup_{t \in I} \{t\} \times \phi_t(a) \subset I \times D^3,
  \]
  and rounding the corners along $\{0, 1\} \times \partial D^3$.
  When the arc $a$ is understood, we simply write $D_{\phi}$ instead of $D_{a, \phi}$.
\end{defn}

It was observed in \cite[Lemma 3.3]{juhász_zemke_dist_slice_disks} that given an (orientation-preserving) self-diffeomorphism $d$ of $(D^3, a)$ such that $d|_{\partial D^3} = \mathrm{id}_{\partial D^3}$, there exists an isotopy $\phi \colon I \times D^3 \to D^3$, such that $\phi_1 = d$. Furthermore, the isotopy class of the deform-spun disk $D_{a, \phi}$ only depends on $d$. Hence we will denote $D_{a,\phi}$ by $D_{a, d}$ for simplicity. Moreover, by \cite[Theorem 5.1]{juhász_zemke_dist_slice_disks}, the element $t_{D_{a, d}}$ is determined by the action of $d_*$ on the knot Floer complex.

\begin{defn}\label{def: deform spun}
Let $K$ be a knot in $S^3$, and suppose that the open 3-ball $B$ intersects $K$ in an unknotted arc.
Then $(S^3 \smallsetminus B, K \smallsetminus B)$ is diffeomorphic to a ball-arc pair $(D^3, a)$.
Suppose that we are given a diffeomorphism $d \in \mathrm{Diff}(S^3, K)$ that is the identity on $B$.
Then the \emph{deform-spun} slice disk $D_{K, d} \subset B^4$ for $-K \# K$ is defined
to be $D_{a, d|_{S^3 \smallsetminus B}}$.
\end{defn}

Now we are able to prove \Cref{cor:1}.

\begin{proof}[Proof of \Cref{cor:1}]
Consider the diffeomorphism $R^\pi:(S^3,K\# K)\rightarrow (S^3,K\# K)$ which swaps the two summands of $K\# K$ as in \Cref{fig:KsharpK}. The action of $R^\pi$ on $\widehat{HFK}(S^3,K\# K)$ (or more generally, $HFK^\infty(S^3,K)$) was identified in \cite[Theorem 8.1]{juhasz_zemke_stabilization_bounds} as
\[
(R^\pi)_\ast = \mathrm{Sw}\circ (1\otimes (1+\Psi\Phi) + \Psi\otimes \Phi),
\]
where $\mathrm{Sw}$ denotes the automorphism of the group
\[
\widehat{HFK}(S^3,K\# K)\simeq \widehat{HFK}(S^3,K)\otimes \widehat{HF}(S^3,K)
\]
which switches the two $\widehat{HFK}(S^3,K)$ factors.

\begin{figure}[h]
    \centering
    \includegraphics[scale=.4]{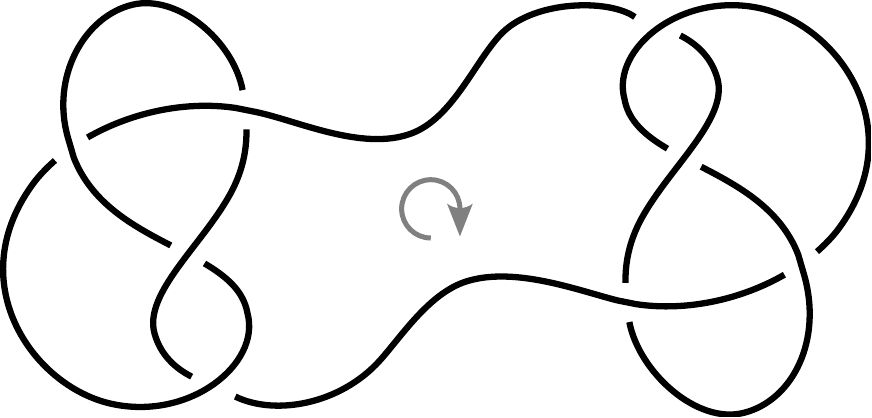}
    \caption{The knot $K\# K$, when $K$ is the figure-eight knot in this figure. The diffeomorphism $R^\pi$ acts by 180$^\circ$ rotation in the plane.} %along the dotted axis.}
    \label{fig:KsharpK}
\end{figure}

Choose any nonzero element $a_0\in \widehat{HFK}(S^3,K)$. If $\Phi(a_0)=0$, then we define $a=a_0$; otherwise we take $a=\Phi(a_0)$, so that $\Phi(a)=\Phi^2(a_0)=0$ since $\Phi^2=0$ \cite[Lemma 4.9]{zemke_linkcob}. Since the hat version of knot Floer homology detects the unknot \cite{os_genusbounds}, we can find a nonzero element $b\in \widehat{HFK}(S^3,K)$ which is different from $a$. Then we have
\[
(R^\pi)_\ast(b\otimes a) = \mathrm{Sw}(b\otimes a) = a\otimes b.
\]
Since $a\ne b$ and $a,b$ are both nonzero, we see that $(R^\pi)_\ast$ is not homotopic to the identity. From \cite[Theorem 5.1]{juhász_zemke_dist_slice_disks}, we then see that the deform-spun slice disks $D_{K\# K,R^\pi}$ and $D_{K\# K,\mathrm{id}}$ of $K\# K\# -K\# -K$, induced by the diffeomorphisms $R^\pi$ and $\mathrm{id}$, satisfy $t_{D_{K\# K,R^\pi}} \ne t_{D_{K\# K,\mathrm{id}}}$ in $\widehat{HFK}(S^3,K\# K\# -K\# -K)$. Therefore, by \Cref{thm:1}, the slice disks $\Lev_{J,s}(D_{K\# K,R^\pi})$ and $\Lev_{J,s}(D_{K\# K,\mathrm{id}})$ are topologically isotopic but not smoothly isotopic.
\end{proof}

\subsection{Cables and Stabilization Distance}\label{sec: cables}

We will now use the same strategy to prove \Cref{thm:1-1}. The type A structure for the $(p, 1)$ cabling pattern, $C_{p,1}$, was computed in \cite{homcabling} (here, we consider the case $p\ge 1$). $CFA^-(S^1 \times D^2,C_{p,1})$ is generated by $a,b_1,\cdots,b_{2p-2}$, and the $A_\infty$ operations are given as follows:
\[
\begin{split}
    m_{3+i}(a,\rho_3,\overbrace{\rho_{23},\cdots,\rho_{23}}^{i},\rho_2) = U^{pi+p} a, &\quad i\ge 0, \\
    m_{4+i+j}(a,\rho_3,\overbrace{\rho_{23},\cdots,\rho_{23}}^{i},\rho_2,\overbrace{\rho_{12},\cdots,\rho_{12}}^{j},\rho_1) = U^{pi+j+1} b_{j+1}, &\quad 0\le j\le p-2,\,i\ge 0, \\
    m_{2+j} (a,\overbrace{\rho_{12},\cdots,\rho_{12}}^{j},\rho_1) = b_{2p-j-2}, &\quad 0 \le j \le p-2, \\
    m_1 (b_j) = U^{p-j} b_{2p-j-1}, &\quad 1 \le j \le p-1, \\
    m_{3+i} (b_j,\rho_2,\overbrace{\rho_{12},\cdots,\rho_{12}}^{i},\rho_1) = U^{i+1}b_{j+i+1}, &\quad 1 \le j \le p-2, \, 0 \le i\le p-j-2, \\
    m_{3+i} (b_j,\rho_2,\overbrace{\rho_{12},\cdots,\rho_{12}}^{i},\rho_1) = b_{j-i-1}, &\quad p+1 \le j \le 2p-2, \, 0 \le i \le j-p-1.
\end{split}
\]
Theorem \ref{thm:1-1} makes uses of the $\HFK^-$ version of the pairing theorem, which states that
\[
CFK^-(S^3,K) \simeq CFA^-(S^1 \times D^2,P) \boxtimes \widehat{CFD}(S^3 \setminus K),
\]
for a pattern $P$ in the solid torus.

\begin{proof}[Proof of \Cref{thm:1-1}]
By \Cref{Theorem: satellite formula for concordance maps}, we have two type D morphisms
\[
f_1,f_2:\widehat{CFD}(S^3 \setminus U)\rightarrow \widehat{CFD}(S^3 \setminus K)
\]
which compute the maps associated to satellites of $D_1$ and $D_2$. Let $x=f_1(v)+f_2(v)$, where $v$ is the single generator of $\widehat{CFD}(S^3 \setminus U)$, and write 
\[
x = \theta + \sum_{J\in \{1,2,3,12,23,123\}} \rho_J \theta_J,
\]
as in the proof of \Cref{lem: no cancellation}. Since $t_{D_1}+t_{D_2}\ne 0$ in $\widehat{HFK}(S^3,K)$, we know that $\theta$ is nonzero. Since the homology of the chain complex
\[
\mathbb{F}_2[U] \simeq CFK^-(S^3,U)\simeq CFA^-(S^1 \times D^2,C_{p,1})\boxtimes \widehat{CFD}(S^3 \setminus U)
\]
is generated by $a\otimes v$, we know that $t_{(D_1)_{p,1}} + t_{(D_2)_{p,1}}$ is the homology class of
\[
(\mathrm{id}_{CFA^-(S^1 \times D^2,C_{p,1})}\boxtimes (f_1+f_2))(a \otimes v),
\]
which contains the term $a \otimes \theta$. Since any $A_\infty$ operation in $CFA^-(S^1 \times D^2,C_{p,1})$ which has some $\mathbb{F}_2[U]$-multiple of $a$ as an output has coefficients of the form $U^{pi+p}$ for some $i\ge 0$, it is clear from the no-cancellation argument of \Cref{lem: no cancellation} that $U^k (t_{(D_1)_{p,1}} + t_{(D_2)_{p,1}})$ is nonvanishing in homology whenever $0 \le k < p$. It then follows from \cite[Theorem 1.1]{juhasz_zemke_stabilization_bounds} that the stabilization distance between $(D_1)_{p,1}$ and $(D_2)_{p,1}$ is at least $p$.
\end{proof}

\section{An infinite family of exotic slice disks} \label{sec:infinite}

As discussed in \S\ref{sec: exotically knotted surfaces}, the techniques developed above are insufficient to distinguish infinite families of pairwise exotic slice disks. Nevertheless, those methods suggest that $\Wh(4_1 \smallsum 4_1)$ should bound  infinitely many exotic slice disks. Instead, we are able to detect this family via an indirect application of Seiberg-Witten invariants.

\begin{figure}[h]
\center
\includegraphics[width=.6\linewidth]{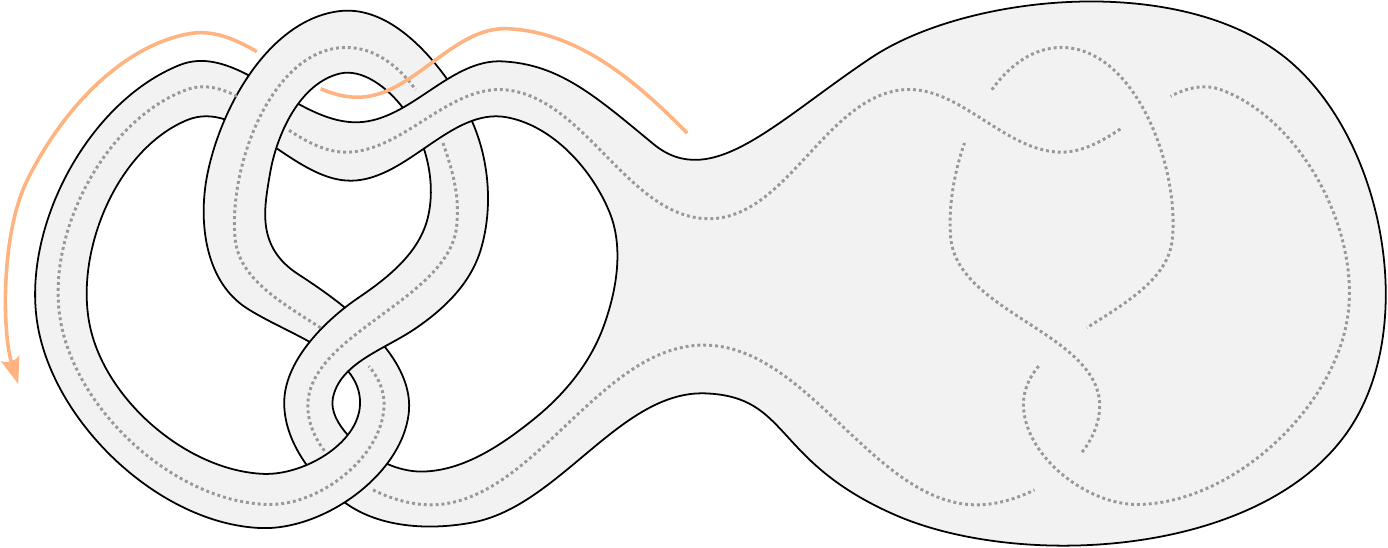}
\caption{An incompressible swallow-follow torus in the complement of $K=4_1 \smallsum -4_1$, decorated with an arrow to indicate the swallow-follow isotopy.}\label{fig:swallow-follow}
\end{figure}

 \begin{proof}[Proof of \Cref{thm:2}]
 Let $K$ denote the knot $4_1 \# -4_1$, and let $D$ denote its standard slice disk.  Consider the isotopy from $K$ to itself given by a longitudinal twist along one of the swallow-follow tori in $S^3 \setminus K$ as depicted in Figure~\ref{fig:swallow-follow}. Dragging $D$ along by any extension of this isotopy to $B^4$ and repeating $n$ times, we obtain a family of ribbon disks $D_n$ for $K$. These disks $D_n$ themselves are not exotic, and they can be distinguished using elementary techniques (e.g., by their inclusion-induced peripheral maps $\pi_1(S^3 \setminus K)\to \pi_1(B^4 \setminus D_n)$, cf \cite{fox:rolling,juhász_zemke_dist_slice_disks}). However, by passing indirectly through a construction due to Gompf \cite{gompf:infinite}, they can also be distinguished using Seiberg-Witten invariants.  We will show that this latter perspective enables us to distinguish the doubled disks $\Wh(D_n)$, which are topologically isotopic rel boundary by \cite{conway_powell_characterisation_homotopy_ribbon}.

\begin{figure}\center
\def\svgwidth{\linewidth}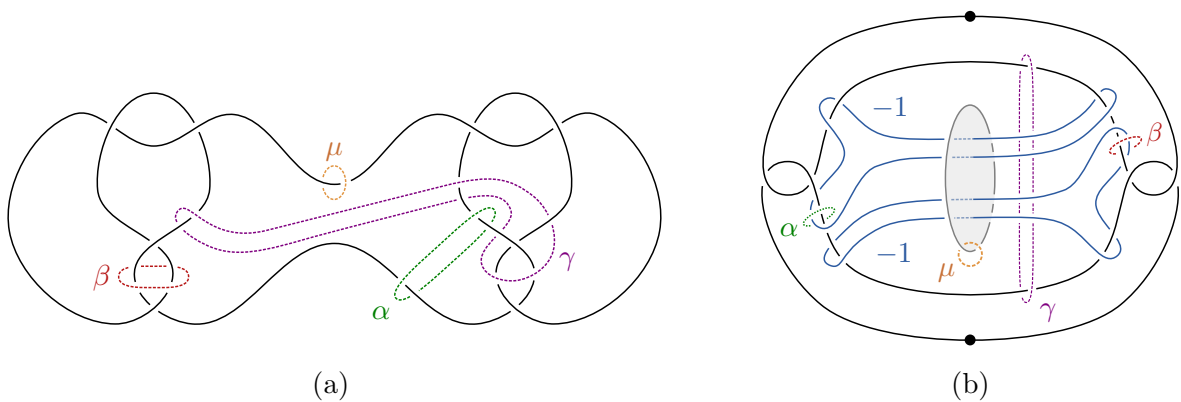
\caption{(a) Decorated curves in the complement of $K=4_1 \smallsum -4_1$. (b) A nonstandard handle diagram for $B^4$ in which $K$ is unknotted and its slice disk (shown in gray) is standard.}\label{fig:akbulut1}
\end{figure}
 
   In \cite{gompf:infinite} (cf \cite{gompf:infinite-handle,akbulut:infinite}), Gompf shows that the contractible 4-manifold $W$ obtained from the exterior of $D \subset B^4$ by attaching a $(-1)$-framed  2-handle along a meridian of $K$ is an infinite-order cork, where the twisting automorphism of $\partial W$ is induced by the aforementioned torus twist in $S^3 \setminus K$. This can be reframed using the disks $D_n$: Given an embedding of $W$ into a larger 4-manifold $Z$, cutting out $W \subset Z$ and regluing by the torus twist is equivalent to simply replacing $W \subset X$ with the 4-manifold $W_n$ given by attaching a $(-1)$-framed meridional 2-handle to the exterior of $D_n$ in $B^4$.

We recall the relevant parts of Gompf's argument, borrowing aspects of Akbulut's exposition in \cite{akbulut:infinite}. Part (a) of Figure~\ref{fig:akbulut1} shows $K$ decorated with four curves $\alpha$, $\beta$, $\gamma$, and $\mu$, the last of which is  a meridian of $K$. Part (b) of Figure~\ref{fig:akbulut1} recasts the situation using a nonstandard handle diagram of $B^4$ in which $K$ and $D$ are unknotted. With framings specified relative to this nonstandard diagram, attach 2-handles to  $B^4 \setminus \nu D$ along $\alpha$ and $\beta$ with framing $-1$ and along $\gamma$ with framing 0. (Note that, in the standard diagram in part (a), $\alpha$ and $\beta$ would be $0$-framed and $\gamma$ would be $(-2)$-framed.) %In addition, we will attach a $(-1)$-framed 2-handle along a copy of $\mu$. However, we leave the curve $\mu$ dashed in the diagram so that we may track it more easily.
It can also be shown that the resulting boundary contains a nonseparating 2-sphere. (This is not obvious from the present diagram but becomes clearer upon using the 2-handles attached along $\alpha$ and $\beta$ to cancel the 1-handles from our nonstandard handle structure on $B^4$. One may then perform simple handleslides to free up a 0-framed 2-handle attached along an unknot that is split from the rest of the diagram.) We may attach a 3-handle along this nonseparating 2-sphere, yielding a new 4-manifold we denote by $X$. At present, since we have not pinned down the precise embedding of the nonseparating 2-sphere, we do not claim that $X$ is well-defined. However, this 2-sphere is disjoint from the meridian curve $\mu$, and we will later attach a 2-handle along $\mu$, making the 4-manifold simply connected. The result of attaching the 3-handle will then become  well-defined by work of Trace \cite{trace}.

% Although we have not pinned down the embedding of the nonseparating 2-sphere, the resulting 4-manifold is well-defined by work of Trace \cite{trace}.

\begin{figure}\center
\def\svgwidth{\linewidth}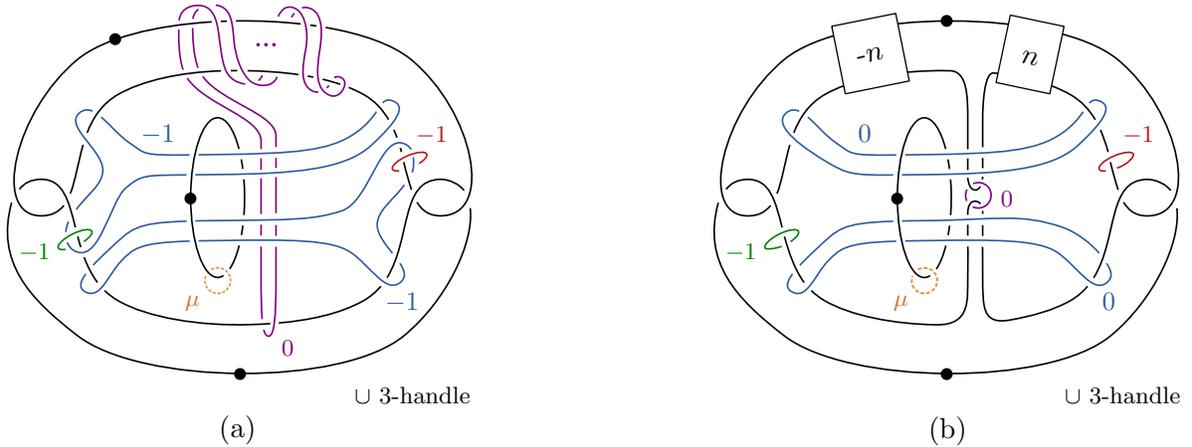
\caption{Diagrams for the 4-manifold $X_n$ obtained by attaching additional handles to the exterior of the disk $D_n \subset B^4$. }\label{fig:akbulut2}
\end{figure}

Let $X_n$ denote the 4-manifold obtained from $X$ by replacing $W \subset X$ with $W_n$. Equivalently, we can apply the torus twist $n$ times to the attaching curves $\alpha$, $\beta$, and $\gamma$ before attaching the 2-handles. Observe that $\gamma$ is the only one of these curves that has essential intersections with the swallow-follow torus, hence we may arrange for the other attaching curves (as well as $\mu$) to be fixed by the torus twist.  The resulting 4-manifold $X_n$ can be shown to have the handle diagram shown in Figure~\ref{fig:akbulut2}(a) (cf \cite[Figure 11]{gompf:infinite-handle}, \cite[Figure 5]{akbulut:infinite}). After further manipulation (namely straightening out $\gamma$ and sliding the original 2-handles from $B^4$ over the 2-handles attached along $\alpha$ and $\beta$), we obtain the diagram in Figure~\ref{fig:akbulut2}(b). %By themselves, the 0-handle, the two 1-handles, and the 0-framed 2-handle attached along $\gamma$ represent the complement of a ribbon disk for $\kappa_n \#-\kappa_n$, where . This allows us to turn Figure~\ref{fig:akbulut2}(b)into the diagram on the left side of Figure~\ref{fig:akbulut-Sn}.  
Using the 0-framed 2-handle attached along $\gamma$, we can combine the two dotted curves representing genuine 1-handles in Figure~\ref{fig:akbulut2}(b) and replace them with a single knotted curve with a dot. This represents the obvious ribbon complement for the slice knot $\kappa_n \# -\kappa_n$, where $\kappa_n$ is the $n$-twist knot. See the left side of Figure~\ref{fig:akbulut-Sn}. This is precisely the result of applying the knot surgery operation of Fintushel-Stern \cite{fintushel-stern:knots} (using the twist knot $\kappa_n$) to the core torus in the simplified picture of $X_0$ shown on the righthand side of Figure~\ref{fig:akbulut-Sn}.

\begin{figure}\center
\def\svgwidth{\linewidth}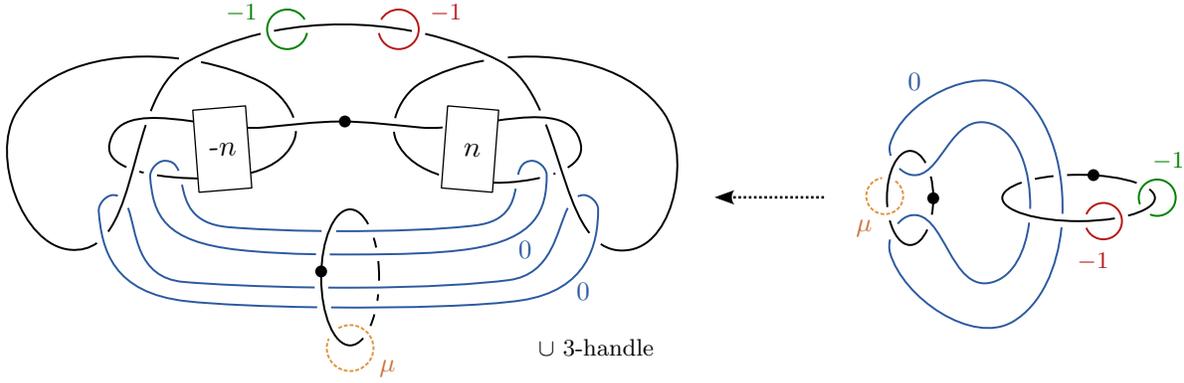
\caption{Realizing $X_n$ as knot surgery on $X_0$.}\label{fig:akbulut-Sn}
\end{figure}

\begin{figure}\center
\def\svgwidth{\linewidth}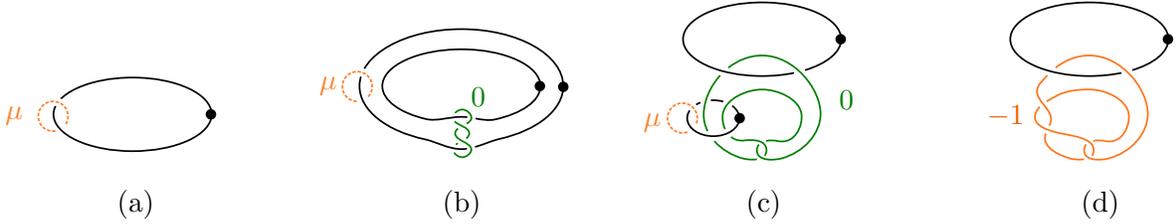
\caption{Part (a) represents the exterior of the unknot's standard slice disk in $B^4$; 2-handles are allowed to pass over the associated 1-handle, but are not pictured. To obtain (b), we replace the exterior of the unknot's standard slice disk with the exterior of its Whitehead double. To obtain (c), we slide the outer 1-handle over the inner 1-handle, then simplify by isotopy.  In (d), we illustrate the effect of attaching a $(-1)$-framed 2-handle along the meridian $\mu$ from part (c) and canceling the resulting 1-/2-handle pair.}\label{fig:std-to-wh}
\end{figure}

\begin{figure}\center
\def\svgwidth{\linewidth}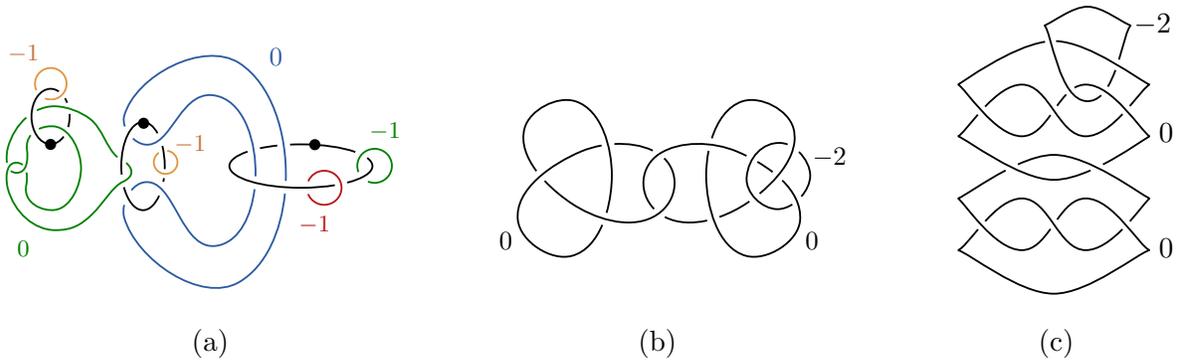
\caption{Diagrams for a 4-manifold obtained by attaching 2-handles to the exterior of $\Wh(D_n) \subset B^4$.}\label{fig:final-mfld}
\end{figure}

We may now turn to the matter at hand: Had we applied this process with the Whitehead doubles $\Wh(D_n)$ instead of $D_n$, we would replace the dotted curve representing $D$ and the meridian curve $\mu$ in the manner shown in parts (a)-(c) of Figure~\ref{fig:std-to-wh}. If we attach a $(-1)$-framed 2-handle along the new meridian, we may simplify to obtain part (d) of Figure~\ref{fig:std-to-wh}. We apply this to the 4-manifolds $X_n$. Moreover, we will also choose to attach a $(-1)$-framed 2-handle along the old meridian. (This realizes the knot surgery torus as a fiber in a cusp neighborhood.) The result of attaching these additional handles to $X_0$ is depicted in part (a) of Figure~\ref{fig:final-mfld}, and the remaining parts (b) and (c) simplify this handle diagram and show that this 4-manifold admits a Stein handle diagram. It follows that this 4-manifold embeds into a (minimal) closed K\"ahler surface $Z$ \cite{lisca-matic:embed}. Since the twist knot $\kappa_n$ has nontrivial Alexander polynomial, the 4-manifolds obtained by performing knot surgery with $\kappa_n$ along the torus in the cusp fiber all have distinct Seiberg-Witten invariants by \cite{fintushel-stern:knots}. (It is clear that the ambient 4-manifold has $b_2^+>1$ and that the torus complement is simply connected, hence we may apply \cite[Theorem 1.5]{fintushel-stern:knots}.) It follows that the slice disks $\Wh(D_n)$ are not smoothly isotopic rel boundary.  \end{proof}

We can upgrade these examples to produce disks that are absolutely exotic, adapting a strategy from \cite{akbulut-ruberman}. 

\begin{proof}[Proof of Corollary~\ref{cor:absolute}]
Let $K'$ denote the knot shown in Figure~\ref{fig:conc-to-wh}.  A calculation in SnapPy \cite{snappy} verifies that $K'$ is hyperbolic with trivial isometry group, which further implies that every diffeomorphism of the pair $(S^3,K')$ is isotopic to the identity (through diffeomorphisms of the pair). Note that $K'$  is obtained from $K=\Wh(4_1\#4_1)$ by a tangle replacement in the center of the diagram. By \cite[Proof of Theorem 2.6]{ruberman}, there exists an invertible concordance $C$ from $K=\Wh(4_1\#4_1)$ to $K'$. Here invertibility means that there is a concordance $C'$ from $K'$ to $K$ such that gluing $C'$ to $C$ along $K'$ yields a product, i.e., $C \cup C'$ and $ K \times [0,1]$ are smoothly isotopic rel boundary. Define a family of slice disks  $D_n'=\Wh(D_n)\cup C$ bounded by $K'$. The slice disks $\Wh(D_n)$ are topologically isotopic rel boundary, hence so are the disks $D'_n$.

\begin{figure}
\includegraphics[width=.7\linewidth]{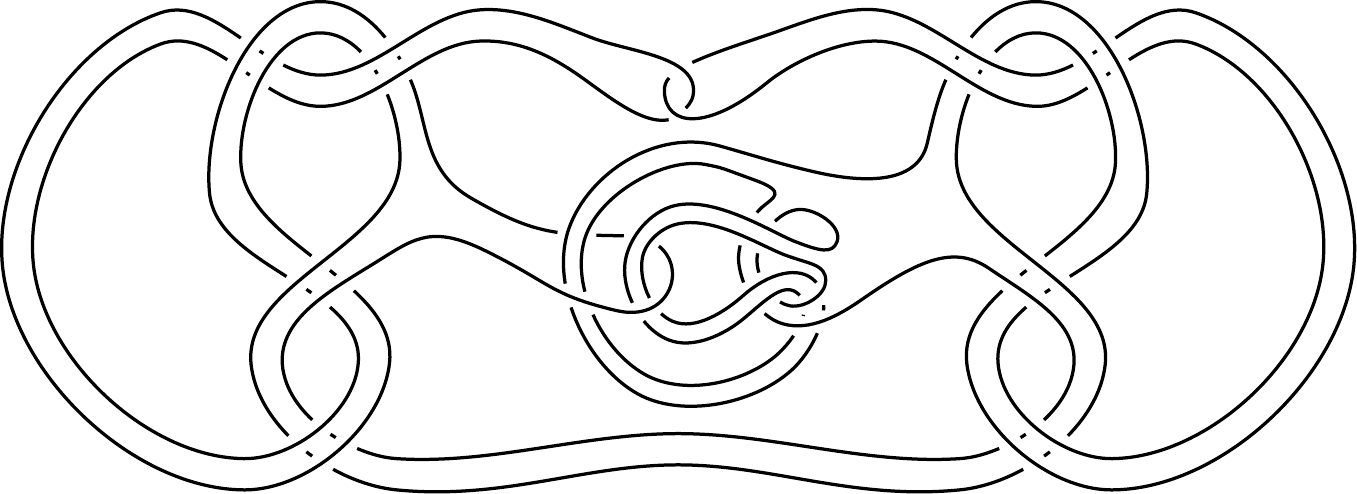}
\caption{A knot $K'$ that is concordant to $\Wh(4_1 \smallsum 4_1)$ and bounds an infinite family of absolutely exotic slice disks.}
\label{fig:conc-to-wh}
\end{figure}

Now suppose that $D'_n$ and $D'_m$ are smoothly isotopic, not necessarily rel boundary, for $n \neq m$. Since $(S^3,K')$ has no nontrivial symmetries, we may further assume that $D'_n$ and $D'_m$ are smoothly isotopic rel boundary. It then follows that $D'_n \cup C'$ and $D'_m \cup C'$ are smoothly isotopic rel boundary. But since $C \cup C'$ is isotopic to the product $K \times [0,1]$ (rel boundary), the disks $D'_n \cup C'$ and $D'_m \cup C'$ can be obtained from $\Wh(D_n)$ and $\Wh(D_m)$ by attaching an extended collar along their boundary. It follows that these extended slice disks are smoothly isotopic rel boundary. However, this implies that the  disks $\Wh(D_n)$ and $\Wh(D_m)$ themselves are smoothly isotopic rel boundary, contradicting the results above.  We conclude that $D_n'$ and $D_m'$ are not smoothly isotopic whenever $n \neq m$.
\end{proof}
\section{Exotic surfaces with diffeomorphic branched  covers}

The goal of this section is to prove Theorem~\ref{thm:dbc}. We begin with a proposition that helps us understand the branched double covers of certain satellite disks. Here we consider a generalization of Levine's patterns $P_{J,k}$ to the case where $k$ is a half-integer. To state it, we recall the definition of \emph{$1/k$-surgery} on a slice disk: Given a smooth slice disk $D \subset B^4$, fix a compact tubular neighborhood $\bar \nu D \cong D \times D^2$ and a choice of meridian $\mu$ of $\partial D \subset S^3$. By taking the disk exterior $B^4 \setminus \nu D$ and attaching a $(-k)$-framed 2-handle along $\mu$,  we
obtain a contractible 4-manifold bounded by $S^3_{1/k}(\partial D)$.  We refer to this 4-manifold
as $1/k$-surgery along the slice disk $D \subset B^4$ and denote it by $B^4_{1/k}(D)$.
%\footnote{\GG{Though for a general $J$, there's probably nothing to do, we can distinguish the covers in the simplest case when $J = U$ and $s > 0.$ \\ }}

\begin{prop}\label{prop:dbc}
\begin{enumerate}[\normalfont\bfseries (a)]
    \item  For any slice disk $D$ in $B^4$, the branched double cover $\Sigma_2(B^4,D_{2,1})$ of the $(2,1)$-cable $D_{2,1}$ of $D$ is diffeomorphic to $B^4_{+1}(D \smallbsum D^r)$.
    \item More generally, let $D \subset B^4$ be a slice disk and $ \mu \in S^3 \setminus \partial D$ a meridian of its boundary. For any knot $J$, half-integer $k \in \tfrac{1}{2}\Z$, and associated doubling pattern $P_{J,k}$, the  branched double cover $\Sigma_2\big(B^4,P_{J,k}(D)\big)$ is diffeomorphic to $\left(B^4 \setminus \nu(D \smallbsum D^r)\right) \cup h$, where $h$ is a $2k$-framed 2-handle attached along $\mu \smallsum J \smallsum J^r$ and $\nu(\, \cdot \,)$ denotes an open tubular neighborhood.
    %[\normalfont\textbf{(b)}]
\end{enumerate}
\end{prop}

\begin{figure}

\

\vspace{.15in}

\def\svgwidth{\linewidth}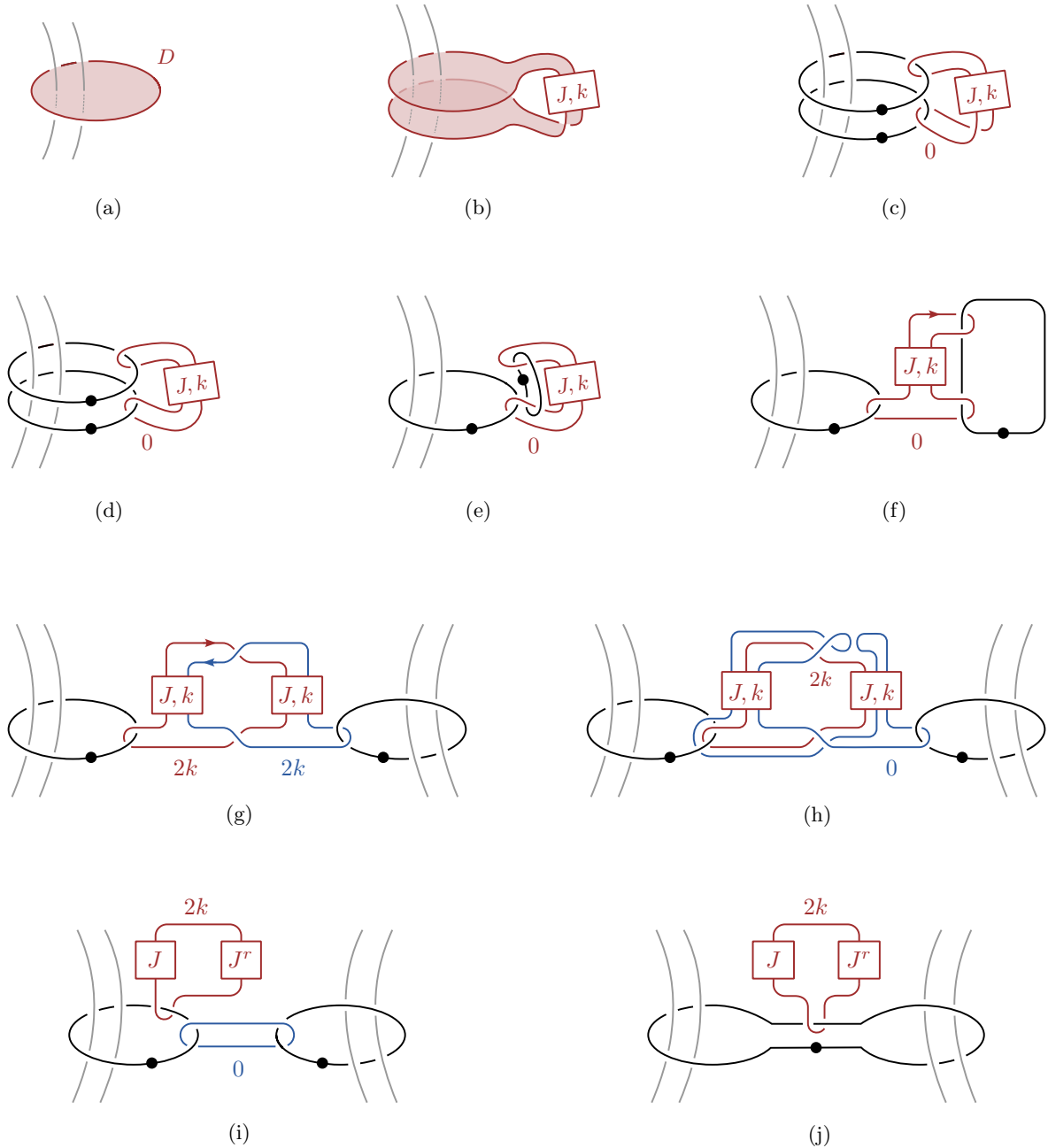

\vspace{.05in}

\caption{(a) The disk $D$ depicted schematically in a  nonstandard handle diagram of $B^4$. (b) The satellite disk $P_{J,k}(D)$. (c)-(f) Handle diagrams for the exterior of $P_{J,k}(D)$ in $B^4$. (g)-(j) Handle diagrams for the branched double cover $\Sigma_2\big(B^4,P_{J,k}(D)\big)$.}
\label{fig:dbc}

\vspace{.15in}
\end{figure}

\begin{proof}
The first claim in the proposition is a special case of the latter claim where $J$ is an unknot and $k=-1/2$, so we focus on the more general claim. The proof is essentially given in Figure~\ref{fig:dbc}, which we unpack below.

Begin by fixing a handle decomposition of $B^4$ in which the original disk $D$ is represented by a standard slice disk bounded by the unknot. This is depicted schematically in part (a) of Figure~\ref{fig:dbc}. The gray arcs represent the other handles that are present in the handle decomposition. In part (b) of Figure~\ref{fig:dbc}, we use this same handle diagram to illustrate the satellite disk $P_{J,k}(D)$.  In part (c), we modify the handle diagram to represent the exterior of $P_{J,k}(D)$ in $B^4$ (see \cite[\S6.2]{GS_4mflds}). The diagram in (d) is obtained by isotopy, (e) is obtained by sliding one 1-handle over the other, and (f) is obtained by further isotopy. In part (g), we depict the branched double cover (see \cite{akbulut-kirby} or \cite[\S6.3]{GS_4mflds}). (We note that we have assumed $2k$ is even in the figure. If $2k$ is odd, then the colors on the bottoms of the red and blue 2-handles are reversed.) The main detail to verify in this step is the framing data for the 2-handles in part (g), which we now explain.

First consider the  framed attaching curve in the center of the diagram in part (f), which we denote by $\gamma$. To pin down what occurs in the $(J,k)$-tangle, fix a diagram of $J$ as a knotted arc with $c$ crossings and writhe $w$. We may construct the $(J,k)$-tangle so that its two strands contribute a total of $4c+2|w|+2|k|$ crossings, where $4c$ crossings come from adding a blackboard-framed pushoff to the underlying diagram of $J$, the next $2|w|$ crossings correct for the difference between the blackboard framing and 0-framing, and  the final $2|k|$ crossings add the desired twisting in the $(J,k)$-tangle. We claim the writhe of $\gamma$ is $2w-2k$: The strands in the tangle have antiparallel orientation, so the first $4c$ crossings cancel in sign; if the strands had parallel orientation, the next $2|w|$ crossings would contribute a signed total of $-2w$ crossings, but instead contribute $2w$ because of the antiparallel orientation; similarly, if the strands had parallel orientations, the final $2|k|$ crossings would contribute $2k$ to the writhe, but instead contribute $-2k$. Thus  the writhe of $\gamma$ is $2w-2k$, hence the 0-framing curve differs from the blackboard framing by $-(2w-2k)=2k-2w$ full signed twists. %the antiparallel orientation would to achieve the desired framing. offset the diagrammatically:  contains $4c+2|w|+2|k|$ crossings, where $c$ is the crossing number of the underlying diagram for $J$ and $w$ is its writhe. As crossings in the attaching curve $\gamma$, the first $4c$ crossings cancel in sign because these strands in $\gamma$ have opposite orientations. If the strands in the $J$-tangle had parallel orientation, then the remaining $2|w|+2|k|$ crossings would amount to a signed count of $-2w$
%\KH{[will fix this awkwardness...]} When the strands in the $J$-tangle are given parallel orientations, there are another $2|w|$ crossings have the opposite sign as $w$. However, in the closed attaching curve $\gamma$, these strands have opposite orientation, so they contribute $2w$ to the writhe of $\gamma$. Similarly, the remaining $|2k|$ crossings would have the same sign as $2k$ if the strands involved had parallel orientations, but their antiparallel orientations result in them contributing $-2k$ to the writhe of the attaching curve $\gamma$. We conclude that the writhe of $\gamma$ is $2w-2k$. It follows that the 0-framing curve differs from the blackboard framing by $2k-2w$ signed full twists.

Now we consider the lift $\tilde \gamma$  of $\gamma$ in the branched double cover, which consists of two framed components. Each component of $\tilde \gamma$ has $2c$ self-crossings (from passing through both $(J,k)$-tangles), with writhe given by $2w$. The framing from $\gamma$ lifts to a framing on each component of $\tilde \gamma$, hence each lifted framing differ from the  blackboard framing by adding $2k-2w$ signed full twists. Since the blackboard framing on each component is $2w$, it follows that the lifted framing on each component of $\tilde \gamma$ is $2k$, as claimed in the diagram in part (g).

Next, we slide one of the components of $\tilde \gamma$ over the other component to obtain the diagram in (h), then perform an isotopy that pulls the modified components out of the $(J,k)$-tangle boxes to obtain the diagram in (i). (To check these steps, it is easiest to work backwards from (i) to (g) by sliding the 0-framed component over the $2k$-framed component so that the former becomes unlinked from the 1-handle on the left.)  %To understand the result of this slide, first note that the components of $\tilde \gamma$ are parallel, hence the modified 2-handle will become unknotted after the slide. If given parallel orientations, the components of $\tilde \gamma$ would have linking number $2k$. Since they are both $2k$-framed, they will become unlinked after the handle slide. Moreover, the modified 2-handle will have framing given by $2k+2k-2(2k)=0$. It is easy to see that it will be unknotted and link the two 1-handle curves as illustrated in part (i) of Figure~\ref{fig:dbc}.
Finally, we slide one of the 1-handles in (i) over the other so that the 0-framed component becomes a 0-framed meridian to whichever dotted circle represents the unmodified 1-handle. We may then slide all other 2-handles off of the latter 1-handle curve until it can be canceled with its 0-framed meridian curve, yielding the diagram in part (j). This final diagram is easily recognized as the result of attaching a $2k$-framed 2-handle to the exterior of $D \smallbsum D^r \subset B^4$ along $\mu \smallsum J \smallsum J^r$, where $\mu$ is a meridian of the disk's boundary.
\end{proof}

We point out that branched covers of satellite disks give rise to many examples of exotic contractible 4-manifolds. This can be true even when the satellite disks themselves are not exotically knotted, as the next example illustrates.
 
\begin{example}\label{ex:positron-cable-dbc} In Example~\ref{ex:946-nonzero}, it was observed that the standard slice disks  $D$ and $D'$ bounded by the knot $K=m(9_{46})$ remain topologically distinct after applying any satellite pattern with nonzero winding number. Nevertheless, for many of the knot doubling patterns $P_{J,k}$, these topologically distinct disks $P_{J,k}(D)$ and $P_{J,k}(D')$ have exotic branched double covers. 

%In Figure~\ref{fig:946-sum}(a), we recall the handle diagram of the exterior of the slice disk $D$ bounded by $K=m(9_{46})$ from Example~\ref{ex:946-nonzero}; reversing the roles of the black and blue curves with a dot-zero swap gives a diagram for the exterior of $D'$. As argued in Example~\ref{ex:946-nonzero}, the curve $\alpha$ bounds a smoothly embedded disk in $B^4 \setminus D'$ but not in $B^4 \setminus D$. 

\begin{figure}[b]
\def\svgwidth{\linewidth}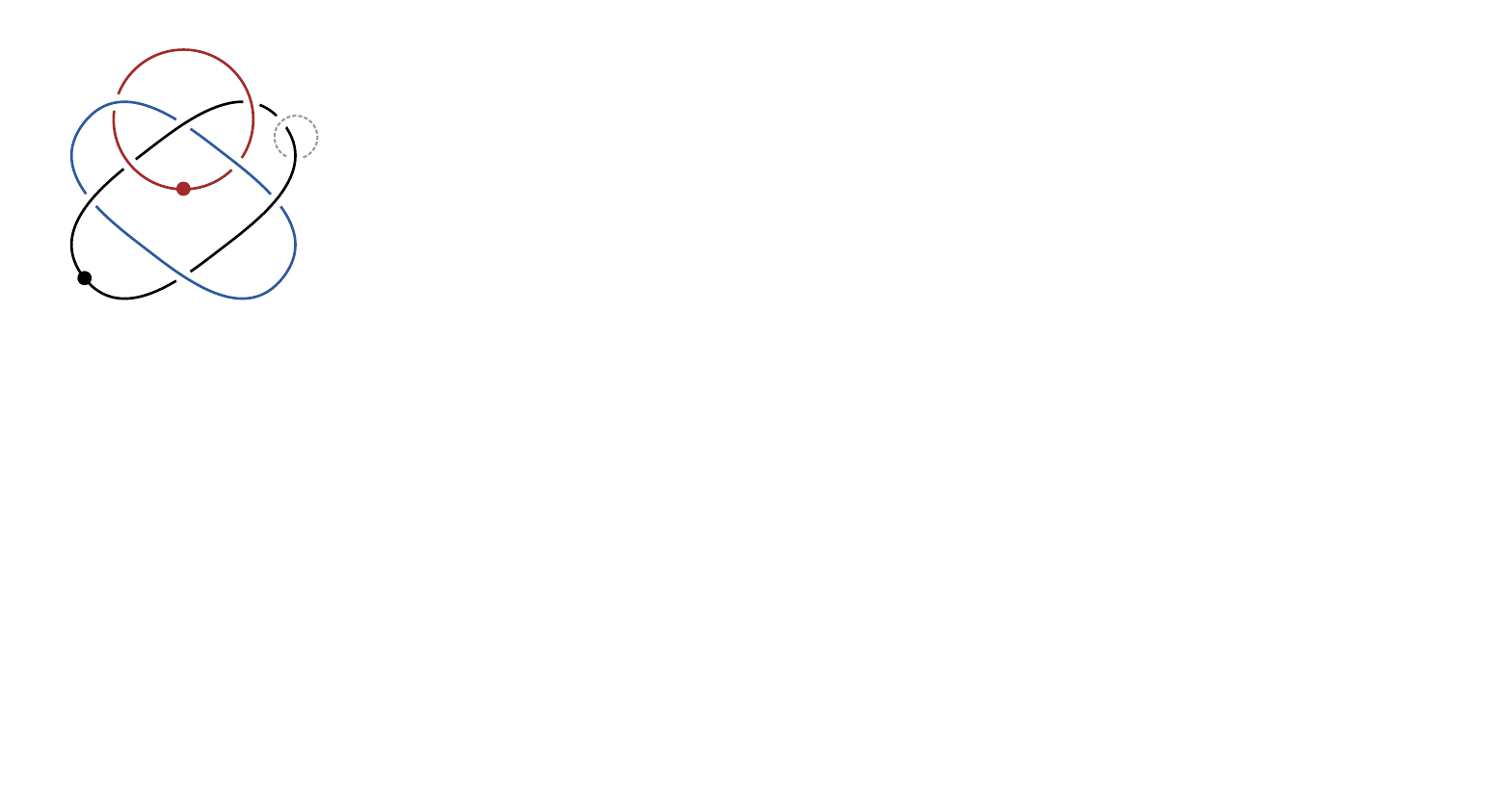
    \caption{}
\label{fig:946-sum}
\end{figure}

To see this, fix a knot $J$ and a half-integer $k \in \tfrac{1}{2} \Z$ such that $J\smallsum J^r$ has maximal Thurston-Bennequin number $\overline{tb}(J\smallsum J^r)>2k$. Following Proposition~\ref{prop:dbc}, we produce a handle diagram for the branched double cover $W=\Sigma_2(B^4,P_{J,k}(D))$ in Figure~\ref{fig:946-sum}(a-c), where $T=J\smallsum J^r$. In particular, the initial diagram in part (a) was obtained earlier in Figure~\ref{fig:946-exterior}, which depicts $D$ as an unknotted disk in a nonstandard the handle diagram of $B^4$. Note that reversing the roles of the black and blue curves with a dot-zero swap gives a diagram for the exterior of $D'$. In turn,  $W'=\Sigma_2(B^4,P_{J,k}(D'))$ is obtained by reversing the roles of the black and blue curves with a dot-zero swap. These branched covers $W$ and $W'$ are contractible 4-manifolds with identical boundary, hence are homeomorphic rel boundary \cite{freedman}.

 To see these 4-manifolds are not diffeomorphic rel boundary, consider attaching $(-1)$-framed 2-handles $h_1$ and $h_2$ to $\partial W=\partial W'$ along the two copies of the curve $\alpha$ in Figure~\ref{fig:946-sum}(c). Since $\alpha$ bounds a smoothly embedded disk in $B^4 \setminus D'$ (as discussed in Example~\ref{ex:946-nonzero}), it is straightforward to see that the enlarged 4-manifold $X'=W'\cup h_1 \cup h_2$ will contain a pair of smoothly embedded 2-spheres of square $-1$. In contrast, consider the other enlarged 4-manifold $X=W \cup h_1 \cup h_2$ shown in Figure~\ref{fig:946-sum}(d). We simplify this by a sequence of handleslides, two handle cancellations, and isotopy to obtain the alternative diagram in part (e). Now choose a Legendrian representative $\mathcal{J} \smallsum \mathcal{J}^r$ of $J \smallsum J^r$ with Thurston-Bennequin number $tb=2k+1$, and let $\mathcal{T}$ be a Legendrian tangle whose closure as in Figure~\ref{fig:946-stein}(a) is $\mathcal{J}\smallsum \mathcal{J}^r$. (If necessary, we can ensure $\mathcal{J}\smallsum \mathcal{J}^r$ has this form by performing Legendrian Reidemeister I moves.) Then $X$ admits a Stein handle structure as depicted in Figure~\ref{fig:946-stein}(b); see \cite{gompf:stein} for background on Stein handle diagrams. Since $X$ admits a Stein structure, it cannot contain any smoothly embedded 2-spheres of square $-1$ \cite{lisca-matic:embed}. Hence we conclude that it is not diffeomorphic to $X'$, and thus $W$ and $W'$ are not diffeomorphic rel boundary.

%\KH{Last mathematical thing to add... just an example showing that, for suitable $J$ and $s$ and where $D$ and $D'$ are the slice disks bounded by $9_{46}$, the branched double covers of $P_{J,s}(D)$ and $P_{J,s}(D')$ are exotic, and can be distinguished using Stein structures / adjunction inequality. Interesting because it's not clear to me that the $(2,1)$-cables of the $9_{46}$ disks are topologically isotopic. Will remark that the same holds for positron disks. This sets up contrast with the proof of Theorem 5. Note that Figure~\ref{fig:positron-cable-dbc} is a placeholder, need to add one more step, then framings.}

\begin{figure}
\def\svgwidth{.85\linewidth}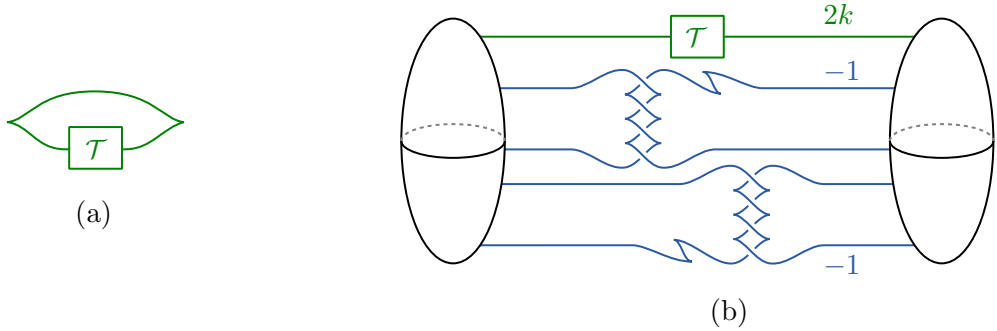
\caption{(a) Representing the Legendrian knot $\mathcal{J}$ as the closure of a Legendrian tangle $\mathcal{T}$. (b) A Stein handle diagram for the 4-manifold $X$.}
\label{fig:946-stein}
\end{figure}

\end{example}

%\begin{remark}
%    \emph{(a)}\ Though the double branched covers of $D_{2,1}$ and $D'_{2,1}$ are exotic, we point out that these cabled disks themselves are not even topologically isotopic. To see this, we first consider the disks $D$ and $D'$ themselves.  If we choose based representatives of the loops $\mu$ and $\beta$ in Figure~\ref{fig:946-sum}, it is straightforward to check that $\pi_1(B^4\setminus D) \cong \langle \mu,\beta : \beta \mu \beta^{-2} \mu^{-1} \rangle$. This admits a surjection onto the symmetric group $S_3$ defined by sending $\mu$ to $(2\ 3)$ and $\beta$ to $(1 \ 2 \ 3 )$ (cf \cite[Proposition 4.4]{akmr:stable}). It follows that $\beta$ is nontrivial in $\pi_1(B^4\setminus D)$, yet it is trivial in $\pi_1(B^4 \setminus D')$ because it bounds a disk in $B^4 \setminus D'$. This implies that $D$ and $D'$ are not  topologically isotopic rel boundary. Moreover, this remains true for the cabled disks $D_{2,1}$ and $D'_{2,1}$. In particular, Remark~\ref{rem:nonzero-winding} points out that there is an inclusion-induced injection $\pi_1(B^4 \setminus \nu D) \hookrightarrow \pi_1(B^4 \setminus D_{2,1})$, hence $\beta$ remains nontrivial in $\pi_1(B^4 \setminus D_{2,1})$ but is trivial in $\pi_1(B^4 \setminus D'_{2,1})$. 
    
    %\emph{(b)} 
%    \KH{[Could potentially point out that the branched covers become diffeomorphic rel boundary after one stabilization using Hayden-Kang-Mukherjee, if we want (where the stabilization is twisted if $2k$ is odd).]}
%\end{remark}

To prove Theorem~\ref{thm:dbc}, we leverage the fact that internal stabilization of knotted surfaces corresponds to external stabilization of their branched covers. In some cases, such as the one given below, exotically knotted surfaces may remain exotic after internal stabilization, yet the exotic branched covers of the original surfaces become diffeomorphic after external stabilization.

\begin{figure}
\def\svgwidth{\linewidth}%% Creator: Inkscape 1.2 (dc2aeda, 2022-05-15), www.inkscape.org
%% Accompanies image file 'positron-and-friends.pdf' (pdf, eps, ps)
\begingroup%
  \makeatletter%
  \providecommand\color[2][]{%
    \errmessage{(Inkscape) Color is used for the text in Inkscape, but the package 'color.sty' is not loaded}%
    \renewcommand\color[2][]{}%
  }%
  \providecommand\transparent[1]{%
    \errmessage{(Inkscape) Transparency is used (non-zero) for the text in Inkscape, but the package 'transparent.sty' is not loaded}%
    \renewcommand\transparent[1]{}%
  }%
  \providecommand\rotatebox[2]{#2}%
  \newcommand*\fsize{\dimexpr\f@size pt\relax}%
  \newcommand*\lineheight[1]{\fontsize{\fsize}{#1\fsize}\selectfont}%
  \ifx\svgwidth\undefined%
    \setlength{\unitlength}{938.53839544bp}%
    \ifx\svgscale\undefined%
      \relax%
    \else%
      \setlength{\unitlength}{\unitlength * \real{\svgscale}}%
    \fi%
  \else%
    \setlength{\unitlength}{\svgwidth}%
  \fi%
  \global\let\svgwidth\undefined%
  \global\let\svgscale\undefined%
  \makeatother%
  \begin{picture}(1,0.25104518)%
    \lineheight{1}%
    \setlength\tabcolsep{0pt}%
    \put(0,0){\includegraphics[width=\unitlength,page=1]{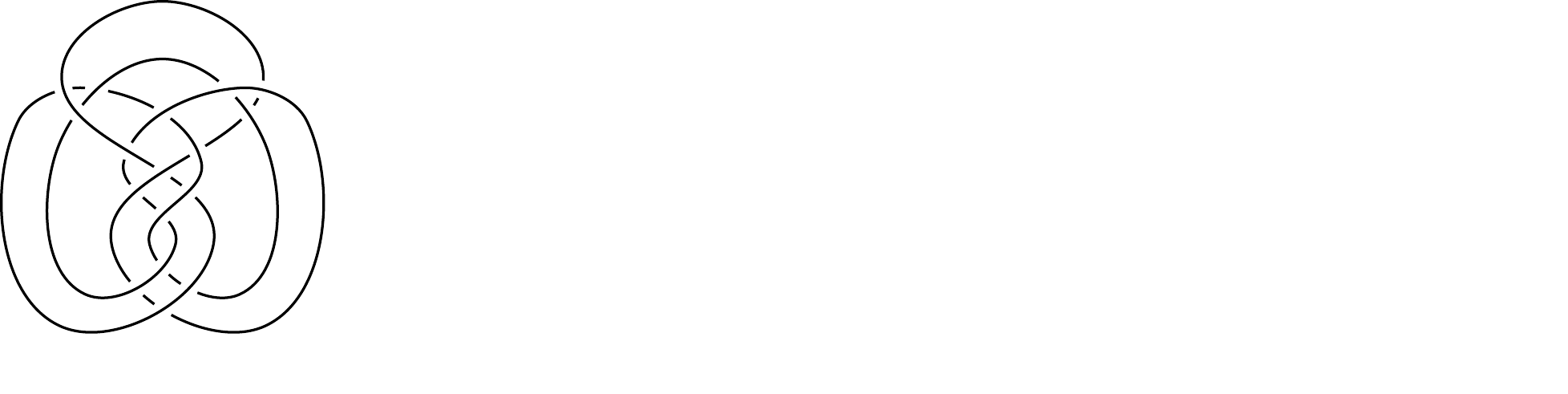}}%
    \put(0.09527329,0.00188333){\color[rgb]{0,0,0}\transparent{0.76999998}\makebox(0,0)[lt]{\smash{\begin{tabular}[t]{l}$K$\end{tabular}}}}%
    \put(0,0){\includegraphics[width=\unitlength,page=2]{positron-and-friends.pdf}}%
    \put(0.35600986,0.00188333){\color[rgb]{0,0,0}\transparent{0.76999998}\makebox(0,0)[lt]{\smash{\begin{tabular}[t]{l}$D$\end{tabular}}}}%
    \put(0,0){\includegraphics[width=\unitlength,page=3]{positron-and-friends.pdf}}%
    \put(0.61646551,0.00188333){\color[rgb]{0,0,0}\transparent{0.76999998}\makebox(0,0)[lt]{\smash{\begin{tabular}[t]{l}$D^r$\end{tabular}}}}%
    \put(0,0){\includegraphics[width=\unitlength,page=4]{positron-and-friends.pdf}}%
    \put(0.88570415,0.00188333){\color[rgb]{0,0,0}\transparent{0.76999998}\makebox(0,0)[lt]{\smash{\begin{tabular}[t]{l}$D'$\end{tabular}}}}%
  \end{picture}%
\endgroup%

\caption{The positron knot $K$ and slice disks that it bounds.}
\label{fig:positron}
\end{figure}

\begin{figure}
\bigskip
\def\svgwidth{\linewidth}%% Creator: Inkscape 1.2 (dc2aeda, 2022-05-15), www.inkscape.org
%% Accompanies image file 'boundary-sum.pdf' (pdf, eps, ps)
\begingroup%
  \makeatletter%
  \providecommand\color[2][]{%
    \errmessage{(Inkscape) Color is used for the text in Inkscape, but the package 'color.sty' is not loaded}%
    \renewcommand\color[2][]{}%
  }%
  \providecommand\transparent[1]{%
    \errmessage{(Inkscape) Transparency is used (non-zero) for the text in Inkscape, but the package 'transparent.sty' is not loaded}%
    \renewcommand\transparent[1]{}%
  }%
  \providecommand\rotatebox[2]{#2}%
  \newcommand*\fsize{\dimexpr\f@size pt\relax}%
  \newcommand*\lineheight[1]{\fontsize{\fsize}{#1\fsize}\selectfont}%
  \ifx\svgwidth\undefined%
    \setlength{\unitlength}{937.78457762bp}%
    \ifx\svgscale\undefined%
      \relax%
    \else%
      \setlength{\unitlength}{\unitlength * \real{\svgscale}}%
    \fi%
  \else%
    \setlength{\unitlength}{\svgwidth}%
  \fi%
  \global\let\svgwidth\undefined%
  \global\let\svgscale\undefined%
  \makeatother%
  \begin{picture}(1,0.23921134)%
    \lineheight{1}%
    \setlength\tabcolsep{0pt}%
    \put(0,0){\includegraphics[width=\unitlength,page=1]{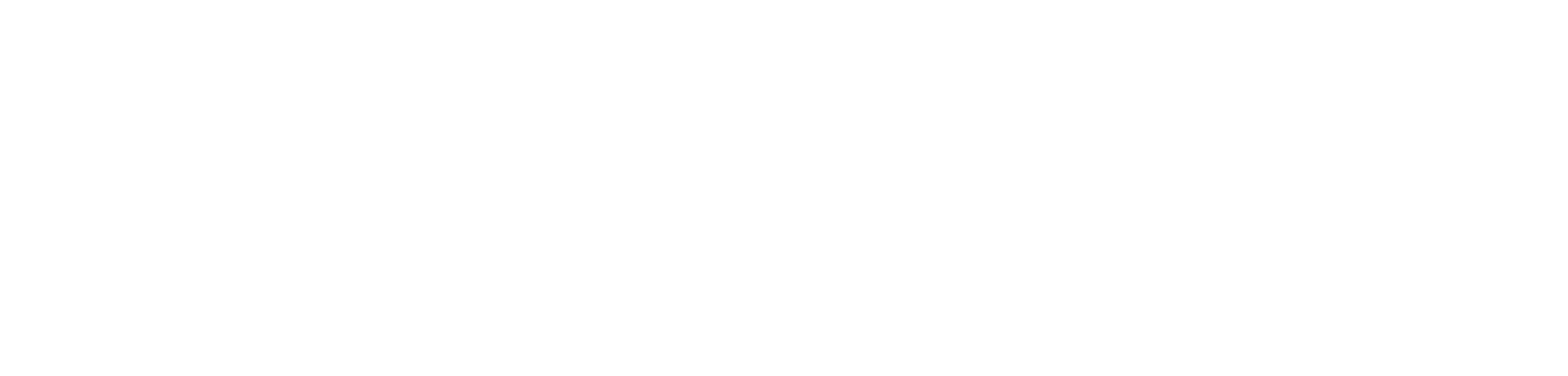}}%
    \put(0.72265143,0.00380091){\color[rgb]{0,0,0}\transparent{0.83421099}\makebox(0,0)[lt]{\smash{\begin{tabular}[t]{l}(b)\end{tabular}}}}%
    \put(0,0){\includegraphics[width=\unitlength,page=2]{boundary-sum.pdf}}%
    \put(0.20260498,0.00386327){\color[rgb]{0,0,0}\transparent{0.83421099}\makebox(0,0)[lt]{\smash{\begin{tabular}[t]{l}(a)\end{tabular}}}}%
    \put(0.51027674,0.22287246){\color[rgb]{0,0,0}\makebox(0,0)[lt]{\smash{\begin{tabular}[t]{l}\raisebox{1pt}{\tiny$+$}{\small $1$}\end{tabular}}}}%
    \put(0,0){\includegraphics[width=\unitlength,page=3]{boundary-sum.pdf}}%
    \put(1.00482442,0.12556758){\color[rgb]{0.50196078,0.50196078,0.50196078}\makebox(0,0)[lt]{\smash{\begin{tabular}[t]{l}{\small $\tau$}\end{tabular}}}}%
    \put(0.7394613,0.23390789){\color[rgb]{0.50196078,0.50196078,0.50196078}\makebox(0,0)[lt]{\smash{\begin{tabular}[t]{l}{\small $\rho$}\end{tabular}}}}%
  \end{picture}%
\endgroup%

\caption{Constructing $D \smallbsum D^r$ and capping it off with the core of the 2-handle in $X=X_{+1}(K \smallsum K^r)$.}
\label{fig:boundary-sum}
\end{figure}

\begin{figure}
\def\svgwidth{\linewidth}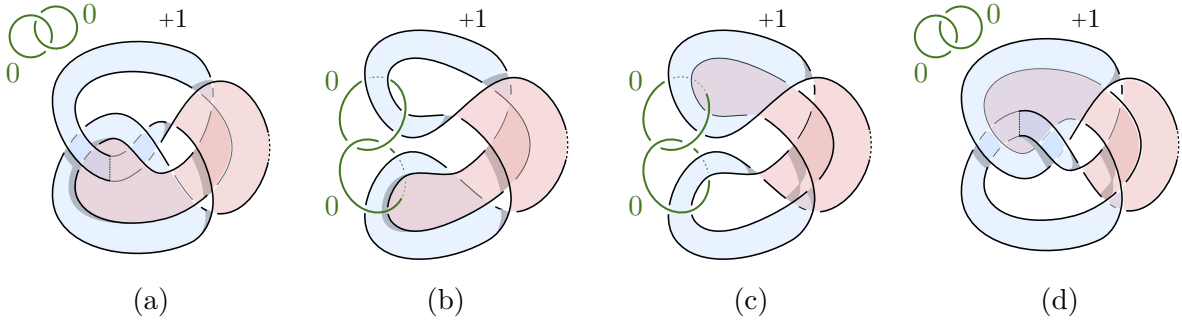
\caption{Applying the key stable isotopy from \cite{akmr:stable}.}
\label{fig:key-isotopy}
\end{figure}

\begin{proof}[Proof of Theorem~\ref{thm:dbc}]
Let $K$ denote the positron knot shown in Figure~\ref{fig:positron}, and let $D$ and $D'$ denote the slice disks shown to its right. By Proposition~\ref{prop:dbc}, the branched double cover $\Sigma_2(B^4,D_{2,1})$ is diffeomorphic to $B^4_{+1}(D \smallbsum D^r)$. The knot $K_{2,1}$ bounds a pair of genus-one surfaces $F,F' \subset B^4$ obtained by internally stabilizing the cabled disks $D_{2,1},D'_{2,1} \subset B^4$, and the surfaces $F$ and $F'$ were shown to be exotically knotted in  \cite{guth_one_not_enough_exotic_surfaces}.

We claim that the branched double covers $\Sigma_2(B^4,F)$ and $\Sigma_2(B^4,F')$ are diffeomorphic rel boundary. To that end, note that the branched double covers of $F=D_{2,1}\smallsum T^2$ and $F'=D'_{2,1}\smallsum T^2$ are diffeomorphic to $\Sigma_2(B^4,D_{2,1}) \smallsum \sxs$ and $\Sigma_2(B^4,D'_{2,1}) \smallsum \sxs$, respectively. By Proposition~\ref{prop:dbc}, the branched covers $\Sigma_2(B^4,D_{2,1})$ and $\Sigma_2(B^4,D'_{2,1})$ are diffeomorphic to $B^4_{+1}(D \smallbsum D^r)$ and $B^4_{+1}(D' \smallbsum (D')^r)$, so we turn our attention to these latter 4-manifolds.

The disk $D \smallbsum D^r$ is illustrated in part (b) Figure~\ref{fig:boundary-sum} (ignoring the label $+1$ for the moment). Note that (up to reversing orientation), the disk $D' \smallbsum (D')^r$ is obtained from $D \smallbsum D^r$ by applying the involution $\tau$ of $B^4$ depicted by $\tau$ in Figure~\ref{fig:boundary-sum}(b). Recall that $B^4_{+1}(D \smallbsum D^r)$ can be constructed as follows: Let  $X=X_{+1}(K \smallsum K^r)$ denote the knot trace obtained by attaching a $(+1)$-framed 2-handle to $B^4$ along $K \smallsum K^r$, and let $S$ and $S'$ denote the smooth 2-spheres in $X$ obtained by capping off $D$ and $D'$, respectively, with the core of the 2-handle in $X$. Then $B^4_{+1}(D \smallbsum D^r)$ and $B^4_{+1}(D' \smallbsum (D')^r)$ are obtained by blowing down $X$ along the $(+1)$-framed 2-spheres $S$ and $S'$, respectively.

To show that $B^4_{+1}(D \smallbsum D^r)$ and $B^4_{+1}(D' \smallbsum (D')^r)$ become diffeomorphic (rel boundary) after a single stabilization with $\sxs$, it suffices to show that $S$ and $S'$ become isotopic (rel boundary) in $X \smallsum \sxs$. We prove this using arguments from \cite{akmr:stable}, with our proof summarized in Figure~\ref{fig:key-isotopy}: In part (a), we see the left half of $D \smallbsum D^r$, viewed as a subset of $S \subset X \smallsum \sxs$. Passing from (a) to (b) uses the ``key stable isotopy'' from \cite[Figure 5]{akmr:stable} (relying on the fact that $S$ has simply connected complement in $X$). The subsurface of $S$ shown in (b) is a disk in $B^4$  bounded by the unknot, hence is isotopic rel boundary to the alternative surface shown in (c). Reversing the ``key  stable isotopy'', we obtain part (d), at which point we have replaced the subset $D \smallbsum D^r$ in $S$ with $D' \smallbsum D^r$. Having freed up the 2-handles corresponding to the $\sxs$-summand of $X \smallsum \sxs$, we can repeat this process and replace the subsurface $D^r$ with $(D')^r$, thus producing an isotopy of $X \smallsum \sxs$ carrying $S$ to $S'$.
\end{proof}

%\begin{remark}
%    As observed in the proof of \Cref{thm:dbc}, under the branched double covering operation, internal stabilizations of smoothly embedded surfaced in $B^4$ translate to stabilization of 4-manifolds. Unlike the case of surfaces, where there are now many examples of exotic pairs that remain exotic after several internal stabilizations,  finding examples of 4-manifolds (with nontrivial boundary) which stay exotic after even a single stabilization is very difficult, and the only known such example was given in \cite{kang2022one}. We point out that exotic 4-manifolds obtained as branched covers of doubled disks often fail \KH{[will come back to finish editing this...]} For doubled disks of the form $\mathrm{Wh}^+(D)$ or $D_{2,1}$, the branched double covers are realized as surgeries along slice disks (per Proposition~\ref{prop:dbc}, hence \cite[Lemma~3.1]{hayden2023one} implies that for any pair $D,D^\prime$ of slice disks bounding a knot $K$, the branched double covers of $\mathrm{Wh}^+(D)$ and $\mathrm{Wh}^+(D^\prime)$ become diffeomorphic after single stabilization with $\sxs$, and the same statement is also true for $D_{2,1}$ and $D^\prime_{2,1}$ under stabilization with $S^2  \tilde \times S^2$.
%\end{remark}

\section{Obstructions from Khovanov homology}

We close by proving Proposition~\ref{prop:khovanov}. For brevity, we assume background at the level of \cite{barnatan}, and we refer the reader to \cite{Hayden_Sundberg_Kh_exotic_disks,Hayden_Kim_Miller_Park_Sundberg_Kh_exotic_seifert,hayden:atomic} for further background on calculations similar to the ones given here.

\begin{proof} 
Let $K=m(9_{46})$ and consider the diagrams of $K$, $K_{2,1}$, and $\Wh(K)$ shown in the top row of Figure~\ref{fig:946-and-doubles}. The bottom row of that figure depicts smoothings of these diagrams,  where 0-resolved crossings are replaced with thin pink arcs (and 1-resolved crossings are unadorned). Note that the smoothings for $K_{2,1}$  and $\Wh(K)$ are essentially obtained by ``doubling'' the smoothing for $K$, except where the additional crossings appear at the top of these diagrams. We obtain Khovanov chain elements by assigning each component of these smoothings the symbol $x$ in the algebra $\mathcal{A}=\F_2\langle 1,x\rangle$. For each of these three chain elements, all 0-resolution arcs join distinct $x$-labeled components, so these chain elements are cycles \cite[Proposition~2.1]{Hayden_Sundberg_Kh_exotic_disks}. 

We claim that the latter two classes are mapped to $1 \in \F_2$ under the cobordism maps $\Kh(D_{2,1})$ and $\Kh(\Wh(K))$. The first steps in these cobordisms are represented by the shaded bands in the diagrams of $K_{2,1}$ and $\Wh(K)$; these indicate band moves which, after performing Reidemeister 1 moves, describe cobordisms from $K_{2,1}$ and $\Wh(K)$ to the link $K_{2,0}$. (Note that these two cobordisms induce different orientations on the link $K_{2,0}$. This will not affect our calculation, so we ignore this detail.) Applying the induced map for these cobordisms (cf \cite[Figure~11]{Hayden_Sundberg_Kh_exotic_disks}), these classes are both carried to the class in $\Kh(K_{2,0})$ shown in the top left of Figure~\ref{fig:movie-1}. From here, these two cobordisms coincide with the doubled disk $D_{2,0}$, and the rest of the cobordism map calculation is summarized in Figures~\ref{fig:movie-1}-\ref{fig:movie-2}. In these figures, we use the label $\Delta$ to represent the split map (cf \cite[Table 1]{Hayden_Sundberg_Kh_exotic_disks}) and the labels $R1$, $R2_a$, and $R2_b$ to denote the chain maps depicted in the first, fifth, and sixth rows of \cite[Table 2]{Hayden_Sundberg_Kh_exotic_disks}, respectively. The final step of each cobordism consists of four local minima capping off the four unknotted components on the top right side of Figure~\ref{fig:movie-2}. Each circle in the associated chain element is $x$-labeled, so the associated cobordism map (cf \cite[Table 1, second row]{Hayden_Sundberg_Kh_exotic_disks}) maps the element to $1 \in \F_2$, as claimed.

\begin{figure}
\includegraphics[width=.9\linewidth]{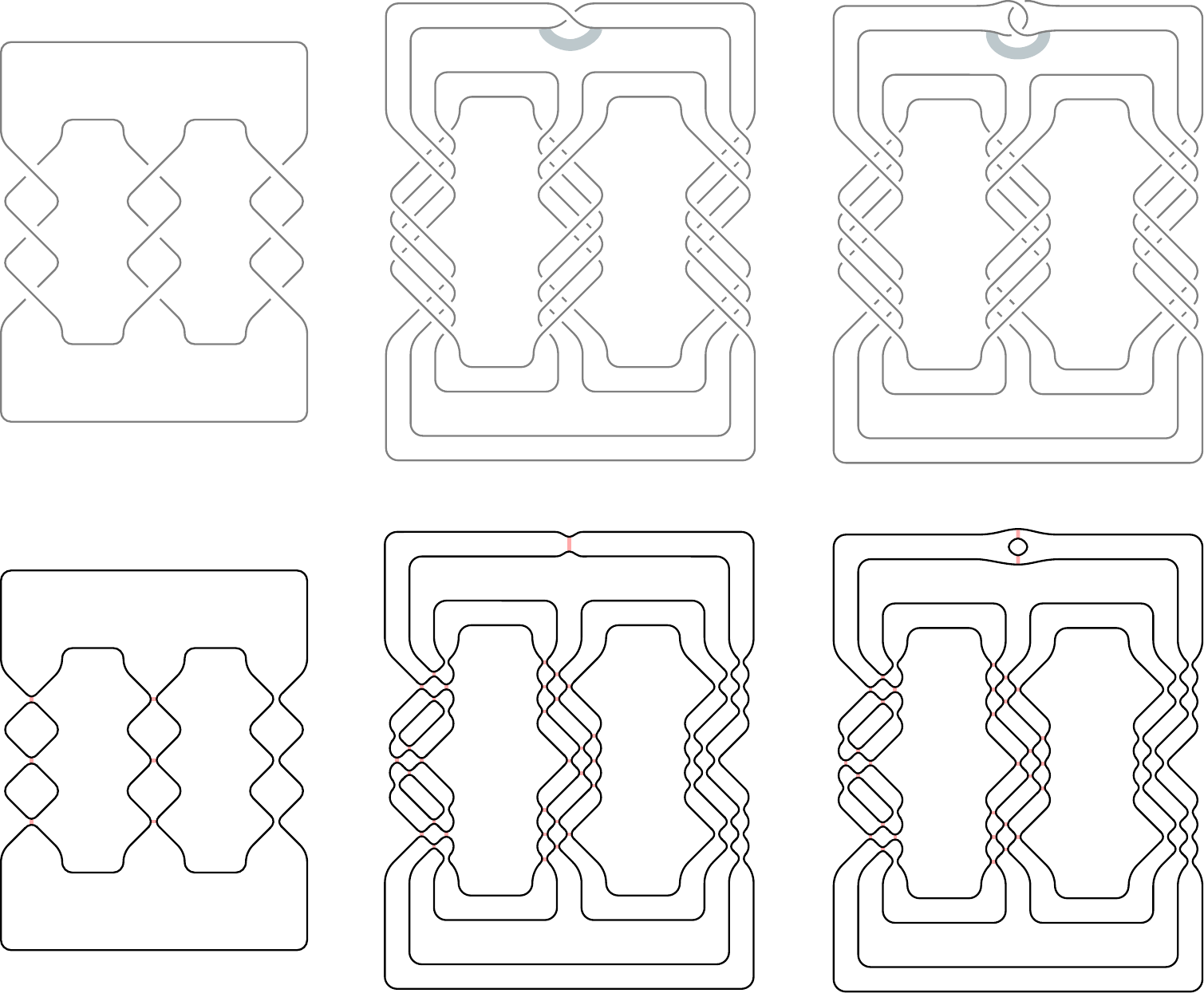}
\caption{The top row presents diagrams for $K=m(9_{46})$, $\Wh(K)$, and $K_{2,1}$. The bottom row depicts certain preferred smoothings of these diagrams.}\label{fig:946-and-doubles}
\end{figure}

\begin{figure}
\def\svgwidth{\linewidth}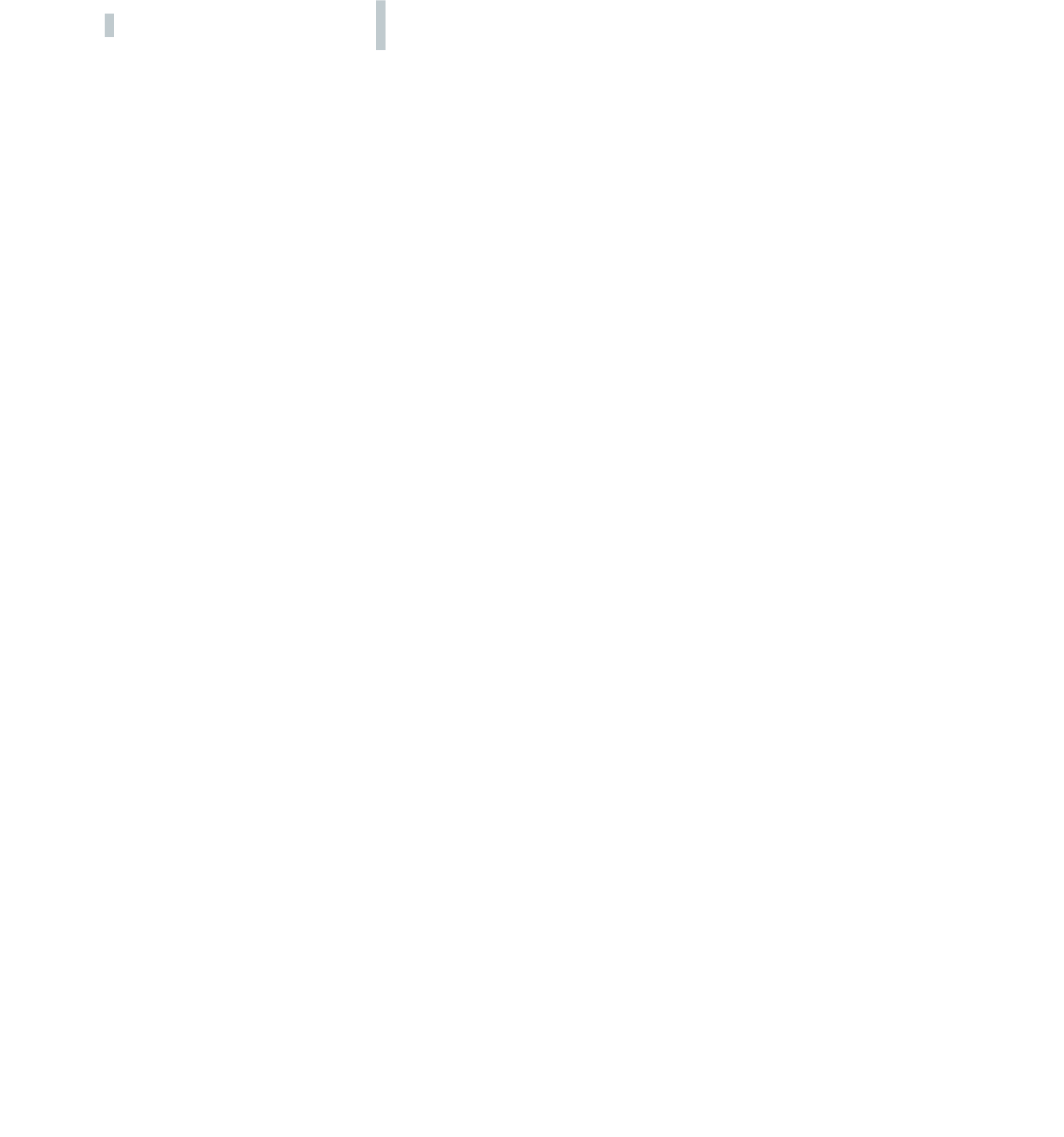
\caption{Calculating the initial steps of the cobordism map.}\label{fig:movie-1}
\end{figure}

\begin{figure}
\def\svgwidth{.7\linewidth}%% Creator: Inkscape 1.2 (dc2aeda, 2022-05-15), www.inkscape.org
%% Accompanies image file '946-double-braid-movie-2.pdf' (pdf, eps, ps)
\begingroup%
  \makeatletter%
  \providecommand\color[2][]{%
    \errmessage{(Inkscape) Color is used for the text in Inkscape, but the package 'color.sty' is not loaded}%
    \renewcommand\color[2][]{}%
  }%
  \providecommand\transparent[1]{%
    \errmessage{(Inkscape) Transparency is used (non-zero) for the text in Inkscape, but the package 'transparent.sty' is not loaded}%
    \renewcommand\transparent[1]{}%
  }%
  \providecommand\rotatebox[2]{#2}%
  \newcommand*\fsize{\dimexpr\f@size pt\relax}%
  \newcommand*\lineheight[1]{\fontsize{\fsize}{#1\fsize}\selectfont}%
  \ifx\svgwidth\undefined%
    \setlength{\unitlength}{704.25088237bp}%
    \ifx\svgscale\undefined%
      \relax%
    \else%
      \setlength{\unitlength}{\unitlength * \real{\svgscale}}%
    \fi%
  \else%
    \setlength{\unitlength}{\svgwidth}%
  \fi%
  \global\let\svgwidth\undefined%
  \global\let\svgscale\undefined%
  \makeatother%
  \begin{picture}(1,0.84556115)%
    \lineheight{1}%
    \setlength\tabcolsep{0pt}%
    \put(0,0){\includegraphics[width=\unitlength,page=1]{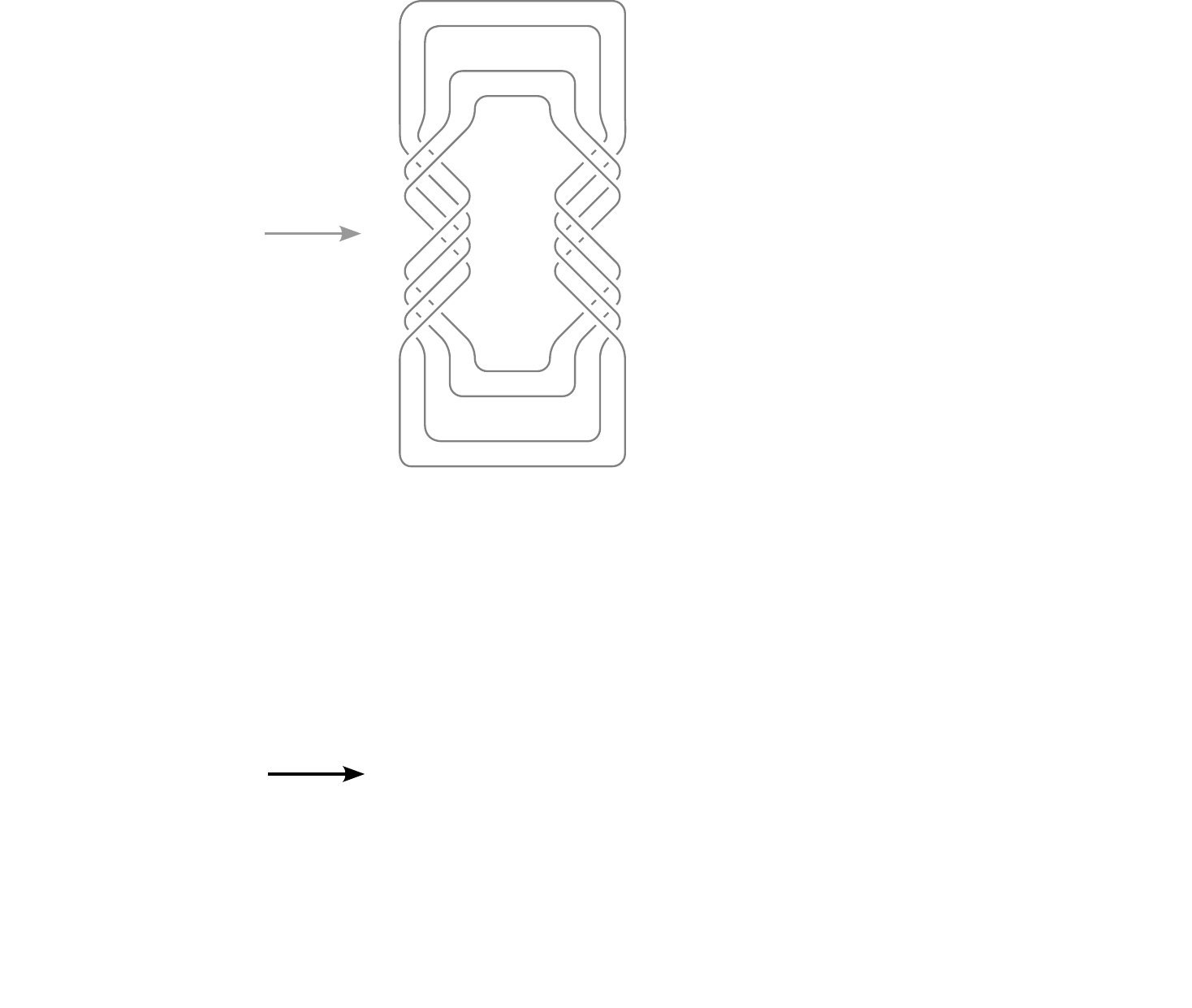}}%
    \put(0.23463568,0.21882869){\color[rgb]{0,0,0}\transparent{0.96078402}\makebox(0,0)[lt]{\smash{\begin{tabular}[t]{l}{\tiny $R2_b$}\end{tabular}}}}%
    \put(0,0){\includegraphics[width=\unitlength,page=2]{946-double-braid-movie-2.pdf}}%
    \put(0.5605174,0.21882869){\color[rgb]{0,0,0}\transparent{0.96078402}\makebox(0,0)[lt]{\smash{\begin{tabular}[t]{l}{\tiny $R2_b$}\end{tabular}}}}%
    \put(0,0){\includegraphics[width=\unitlength,page=3]{946-double-braid-movie-2.pdf}}%
    \put(0.71813215,0.21882869){\color[rgb]{0,0,0}\transparent{0.96078402}\makebox(0,0)[lt]{\smash{\begin{tabular}[t]{l}{\tiny $R2_b$}\end{tabular}}}}%
  \end{picture}%
\endgroup%

\caption{Calculating the final steps of the cobordism map.}\label{fig:movie-2}
\end{figure}

Next we claim that these elements lie in the kernels of the maps $\Kh(D'_{2,1})$ and $\Kh(\Wh(D'))$. These cobordisms begin with the same initial saddle cobordisms to $K_{2,0}$; as discussed above, the induced maps carry the elements in question to the class in $\Kh(K_{2,0})$ depicted in the top left of Figure~\ref{fig:movie-1}. From here, the cobordisms coincide with $D'_{2,0}$. This cobordism continues with the saddle move illustrated on the left side of Figure~\ref{fig:Dprime}. At the chain level, this corresponds to a merge map (cf \cite[Table 1, third row]{Hayden_Sundberg_Kh_exotic_disks}) between two distinct $x$-labeled components of the smoothing, hence vanishes. It follows that the element in question is mapped to $0 \in \F_2$. 

\begin{figure}
\includegraphics[width=.6\linewidth]{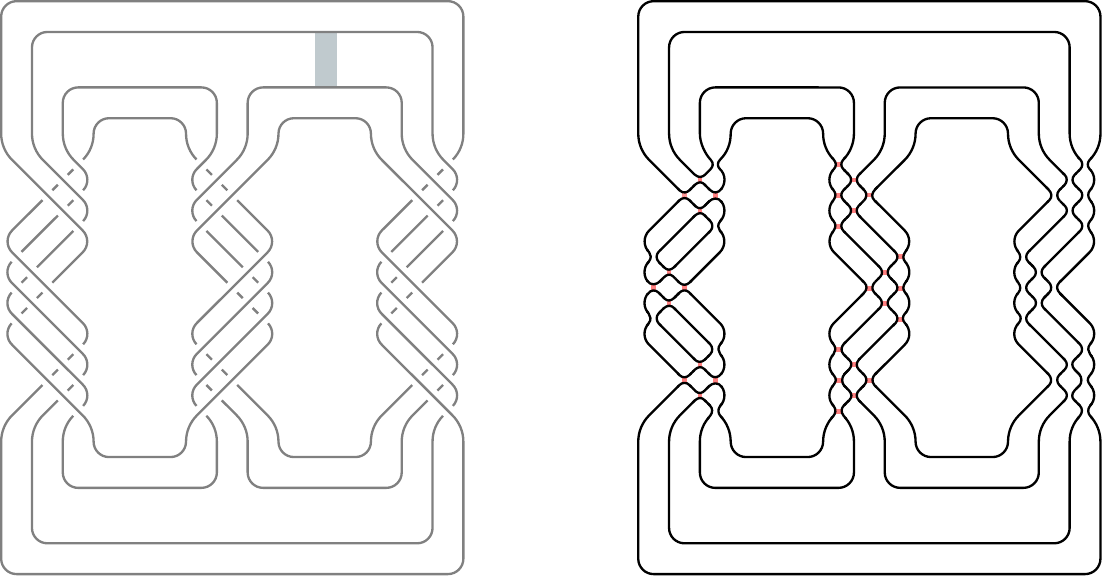}
\caption{The cobordism movie for doubles of $D'$ involve a saddle move along the band shown on the left. This saddle move merges two distinct circles in the smoothing on the right.}\label{fig:Dprime}
\end{figure}

Finally, we claim that the genus-one surfaces obtained by any internal stabilization of $D_{2,1}$ and $D'_{2,1}$ induce distinct maps on (reduced) Bar-Natan homology over $\F_2[U]$. The argument is analogous to the proof of \cite[Proposition~4.2]{hayden:atomic}, so we only outline it: First observe that the calculation establishing $\Kh(D_{2,1}) \neq \Kh(D'_{2,1})$ over $\F_2$ implies $\BNr(D_{2,1}) \neq \BNr(D'_{2,1})$ over $\F_2[U]$ (cf \cite[Proposition 4.1]{hayden:atomic}). As a technical point, we note that working with the reduced theory $\BNr$ requires puncturing our surfaces so that they are cobordisms between $K_{2,1}$ and the unknot. Next, consider the mirrored, time-reversed disks $-D_{2,1}$ and $-D'_{2,1}$. The duality of the cobordism maps under mirroring implies that the maps $\BNr(-D_{2,1}),\BNr(-D'_{2,1}): \F_2[U] \to \BNr(-K_{2,1})$ are distinct. Note that these maps are determined by the images of the generator $1 \in \F_2[U]$. By Euler characteristic considerations, these maps send $1 \in \F_2[U]$ to elements in bigrading $(h,q)=(0,0)$. As verified by computer calculation (in \cite{knotjob}) and depicted in Table~\ref{table:main-ss}, all nonzero elements in $\BNr(-K_{2,1})$ with  bigrading $(0,0)$ survive multiplication by $U$. It follows that the difference $\BNr(-D_{2,1})(1)-\BNr(-D'_{2,1})(1)$ survives multiplication by $U$. Since  internal stabilization  has the effect of multiplying the cobordism maps by $U$ (\cite[Proposition~6.11]{lipshitz-sarkar:mixed}) the claim follows.
\end{proof}

\begin{table}[h!]
\tiny
\centering
\setlength\extrarowheight{2pt}
\begin{tabular}{|c||cc|>{\centering}m{.024\textwidth}|>{\centering}m{.024\textwidth}|>{\centering}m{.02\textwidth}|>{\centering}m{.02\textwidth}|>{\centering}m{.02\textwidth}|>{\centering}m{.02\textwidth}|c|}
 \arrayrulecolor{black}
\multicolumn{10}{c}{Page 1} \smallskip \\
\cline{1-10}  
\color{black}\backslashbox{\!$q$\!}{\!$h$\!} &  $\hdots$   &\hspace{-8pt} {\color{black}{\vrule}} \hspace{-2pt} \color{black}$-3$ \hspace{-5pt} & \color{black}{$-2$}  & \color{black}$-1$ & \color{black}$0$ & \color{black}$1$ & \color{black}$2$ & \color{black}$3$ & \color{black}$4$  \\
\hhline{=||=========}%\cmidrule{1-10}\morecmidrules\cmidrule{1-10}
$4$     &   &   &   &   &   &   &   &   & $1$ \\ 
\hhline{-||~--------}%\hhline{-||~--------}
$2$     &   &   &   &   &   &  &  1 &  1 &   \\
\hhline{-||~--------}
$0$     &   &   &   &   & \cellcolor{cellgray} 2 & 2  &1  & $ $ &   \\
\hhline{-||~--------}
$-2$     &   &   &   & 1   &   3 & 2 &   &   &   \\
\hhline{-||~--------}
$-4$     &   &   &    4 &  5 & 1 &  &   &   &   \\
\hhline{-||~--------}
$-6$  &   & 5  & 4   & 1 &  &   &   &   &   \\
\hhline{-||~--------}
$-8$     &   & 5   & 1 &   &   &   &   &   &   \\
\hhline{-||~--------}
$-10$    &   &  &  &   &   &   &   &   &   \\
\hhline{-||}
$\vdots$ & \  \reflectbox{$\ddots$}
\\
\end{tabular}\hfill
\begin{tabular}{|c||cc|>{\centering}m{.024\textwidth}|>{\centering}m{.024\textwidth}|>{\centering}m{.02\textwidth}|>{\centering}m{.02\textwidth}|>{\centering}m{.02\textwidth}|>{\centering}m{.02\textwidth}|c|}
\multicolumn{10}{c}{}\\
\multicolumn{10}{c}{Page 2} \smallskip \\
\cline{1-10}  
\color{black}\backslashbox{\!$q$\!}{\!$h$\!} &  $\hdots$   &\hspace{-8pt} {\color{black}{\vrule}} \hspace{-2pt} \color{black}$-3$ \hspace{-5pt} & \color{black}$-2$ & \color{black}$-1$ & \color{black}$0$ & \color{black}$1$ & \color{black}$2$ & \color{black}$3$ & \color{black}$4$ \\
\hhline{=||=========}
$4$     &   &   &   &   &   &   &   &   & \\
\hhline{-||~--------}
$2$     &   &   &   &   &   &   &   &   &   \\
\hhline{-||~--------}
$0$     &   &   &   &   &   \cellcolor{cellgray} 2  &   &   &   &   \\
\hhline{-||~--------}
$-2$     &   &   &   &   &   &   &   &   &   \\
\hhline{-||~--------}
$-4$     &   &   &   &   1 &   &   &   &   &   \\
\hhline{-||~--------}
$-6$     &   &   &   &  &   &   &   &   &   \\
\hhline{-||~--------}
$-8$     &   &   &   &   &   &   &   &   &   \\
\hhline{-||~--------}
$-10$     &   &   &   &   &   &   &   &   &   \\
\hhline{-||}
$\vdots$ & \ \reflectbox{$\ddots$} \\
\multicolumn{10}{c}{}\\
\end{tabular}

\vspace{-10pt}

\captionsetup{width=.95\linewidth}
\caption{The first two pages of the reduced Bar-Natan--Lee--Turner spectral sequence for $(9_{46})_{2,1}$, shown for  $h \geq -3$ and $q \geq -10$.} 
\label{table:main-ss}
%\vspace{-7.5pt}
\end{table}

\bibliographystyle{amsalpha} 
\bibliography{main}

\newcommand{\etalchar}[1]{$^{#1}$}
\providecommand{\bysame}{\leavevmode\hbox to3em{\hrulefill}\thinspace}
\providecommand{\MR}{\relax\ifhmode\unskip\space\fi MR }
\providecommand{\MRhref}[2]{%
  \href{http://www.ams.org/mathscinet-getitem?mr=#1}{#2}
}
\providecommand{\href}[2]{#2}
\begin{thebibliography}{HKM{\etalchar{+}}22}

\bibitem[AK80]{akbulut-kirby}
Selman Akbulut and Robion Kirby, \emph{Branched covers of surfaces in
  {$4$}-manifolds}, Math. Ann. \textbf{252} (1979/80), no.~2, 111--131.
  \MR{593626}

\bibitem[Akb17]{akbulut:infinite}
Selman Akbulut, \emph{On infinite order corks}, Proceedings of the {G}\"{o}kova
  {G}eometry-{T}opology {C}onference 2016, G\"{o}kova Geometry/Topology
  Conference (GGT), G\"{o}kova, 2017, pp.~151--157. \MR{3676087}

\bibitem[AKMR15]{akmr:stable}
Dave Auckly, Hee~Jung Kim, Paul Melvin, and Daniel Ruberman, \emph{Stable
  isotopy in four dimensions}, J. Lond. Math. Soc. (2) \textbf{91} (2015),
  no.~2, 439--463. \MR{3355110}

\bibitem[AR16]{akbulut-ruberman}
Selman Akbulut and Daniel Ruberman, \emph{Absolutely exotic compact
  4-manifolds}, Comment. Math. Helv. \textbf{91} (2016), no.~1, 1--19.
  \MR{3471934}

\bibitem[BN05]{barnatan}
Dror Bar-Natan, \emph{Khovanov's homology for tangles and cobordisms}, Geom.
  Topol. \textbf{9} (2005), 1443--1499. \MR{2174270}

\bibitem[BS16]{baykur_sunukjian_knotted_surf_stab}
R.~\.{I}nan\c{c} Baykur and Nathan Sunukjian, \emph{Knotted surfaces in
  4-manifolds and stabilizations}, J. Topol. \textbf{9} (2016), no.~1,
  215--231. \MR{3465848}

\bibitem[CDGW]{snappy}
Marc Culler, Nathan~M. Dunfield, Matthias Goerner, and Jeffrey~R. Weeks,
  \emph{Snap{P}y, a computer program for studying the geometry and topology of
  $3$-manifolds}, \url{http://snappy.computop.org}.

\bibitem[CP21a]{conway_powell_characterisation_homotopy_ribbon}
Anthony Conway and Mark Powell, \emph{Characterisation of homotopy ribbon
  discs}, Adv. Math. \textbf{391} (2021), Paper No. 107960, 29. \MR{4300918}

\bibitem[CP21b]{conway2021embedded}
\bysame, \emph{Embedded surfaces with infinite cyclic knot group}, Geometry and
  Topology (2021).

\bibitem[DHM21]{DaiHeddenMallick_corks_invol_floer}
Irving Dai, Matthew Hedden, and Abhishek Mallick, \emph{Corks, involutions, and
  {H}eegaard {F}loer homology}, arXiv:2002.02326 (2021).

\bibitem[FKV87]{finashin1987exotic}
SM~Finashin, M~Kreck, and OY~Viro, \emph{Exotic knottings of surfaces in the
  4-sphere}, Bulletin of the American Mathematical Society (New Series)
  \textbf{17} (1987), no.~2, 287--290.

\bibitem[Fox66]{fox:rolling}
R.~H. Fox, \emph{Rolling}, Bull. Amer. Math. Soc. \textbf{72} (1966), 162--164.
  \MR{184221}

\bibitem[FQ90]{freedman_quinn_top_4mflds}
Michael~H. Freedman and Frank Quinn, \emph{Topology of 4-manifolds}, Princeton
  Mathematical Series, vol.~39, Princeton University Press, Princeton, NJ,
  1990. \MR{1201584}

\bibitem[Fre82]{freedman}
Michael~Hartley Freedman, \emph{The topology of four-dimensional manifolds}, J.
  Differential Geom. \textbf{17} (1982), no.~3, 357--453. \MR{679066}

\bibitem[FS97]{fintushel1997surfaces}
Ronald Fintushel and Ronald~J Stern, \emph{Surfaces in 4-manifolds},
  Mathematical Research Letters \textbf{4} (1997), no.~6, 907--914.

\bibitem[FS98]{fintushel-stern:knots}
Ronald Fintushel and Ronald~J. Stern, \emph{Knots, links, and {$4$}-manifolds},
  Invent. Math. \textbf{134} (1998), no.~2, 363--400. \MR{1650308}

\bibitem[Glu62]{gluck1962embedding}
Herman Gluck, \emph{The embedding of two-spheres in the four-sphere},
  Transactions of the American Mathematical Society \textbf{104} (1962), no.~2,
  308--333.

\bibitem[Gom98]{gompf:stein}
R.~Gompf, \emph{Handlebody construction of {S}tein surfaces}, Ann. of Math. (2)
  \textbf{148} (1998), no.~2, 619--693.

\bibitem[Gom17a]{gompf:infinite}
Robert~E. Gompf, \emph{Infinite order corks}, Geom. Topol. \textbf{21} (2017),
  no.~4, 2475--2484. \MR{3654114}

\bibitem[Gom17b]{gompf:infinite-handle}
\bysame, \emph{Infinite order corks via handle diagrams}, Algebr. Geom. Topol.
  \textbf{17} (2017), no.~5, 2863--2891. \MR{3704246}

\bibitem[GS99]{GS_4mflds}
Robert~E. Gompf and Andr\'{a}s~I. Stipsicz, \emph{{$4$}-manifolds and {K}irby
  calculus}, Graduate Studies in Mathematics, vol.~20, American Mathematical
  Society, Providence, RI, 1999. \MR{1707327}

\bibitem[Gut22]{guth_one_not_enough_exotic_surfaces}
Gary Guth, \emph{For exotic surfaces with boundary, one stabilization is not
  enough}, arXiv:2207.11847 (2022).

\bibitem[Hay23]{hayden:atomic}
Kyle Hayden, \emph{An atomic approach to {W}all-type stabilization problems},
  arxiv:2302.10127 (2023).

\bibitem[Hed07]{hedden_whitehead_doubles}
Matthew Hedden, \emph{Knot floer homology of whitehead doubles}, Geometry \&
  Topology \textbf{11} (2007), no.~4, 2277 -- 2338.

\bibitem[HKK{\etalchar{+}}21]{hkkmps}
Kyle Hayden, Alexandra Kjuchukova, Siddhi Krishna, Maggie Miller, Mark Powell,
  and Nathan Sunukjian, \emph{Brunnian exotic surface links in the 4-ball},
  Michigan Math. J. \textbf{arXiv:2106.13776} (2021).

\bibitem[HKM{\etalchar{+}}22]{Hayden_Kim_Miller_Park_Sundberg_Kh_exotic_seifert}
Kyle Hayden, Seungwon Kim, Maggie Miller, JungHwan Park, and Isaac Sundberg,
  \emph{Seifert surfaces in the 4-ball}, arXiv:2205.15283 (2022).

\bibitem[Hom14]{homcabling}
Jennifer Hom, \emph{Bordered {H}eegaard {F}loer homology and the tau-invariant
  of cable knots}, J. Topol. \textbf{7} (2014), no.~2, 287--326. \MR{3217622}

\bibitem[HS21]{Hayden_Sundberg_Kh_exotic_disks}
Kyle Hayden and Isaac Sundberg, \emph{Khovanov homology and exotic surfaces in
  the 4-ball}, arXiv:2108.04810 (2021).

\bibitem[JM16]{juhasz2016concordance}
Andr\'{a}s Juh\'{a}sz and Marco Marengon, \emph{Concordance maps in knot
  {F}loer homology}, Geom. Topol. \textbf{20} (2016), no.~6, 3623--3673.
  \MR{3590358}

\bibitem[JZ20]{juhász_zemke_dist_slice_disks}
Andr\'{a}s Juh\'{a}sz and Ian Zemke, \emph{Distinguishing slice disks using
  knot {F}loer homology}, Selecta Math. (N.S.) \textbf{26} (2020), no.~1, Paper
  No. 5, 18. \MR{4045151}

\bibitem[JZ21]{juhasz_zemke_stabilization_bounds}
Andr\'as Juh\'asz and Ian Zemke, \emph{Stabilization distance bounds from link
  {F}loer homology}, arXiv:1810.09158 (2021).

\bibitem[Kir78]{kirby_problem_list}
Rob Kirby, \emph{Problems in low dimensional manifold theory}, Algebraic and
  geometric topology ({P}roc. {S}ympos. {P}ure {M}ath., {S}tanford {U}niv.,
  {S}tanford, {C}alif., 1976), {P}art 2, Proc. Sympos. Pure Math., XXXII, Amer.
  Math. Soc., Providence, R.I., 1978, pp.~273--312. \MR{520548}

\bibitem[Lev12]{levine_doubling_operators}
Adam~Simon Levine, \emph{Knot doubling operators and bordered {H}eegaard
  {F}loer homology}, J. Topol. \textbf{5} (2012), no.~3, 651--712. \MR{2971610}

\bibitem[LM97]{lisca-matic:embed}
P.~Lisca and G.~Mati\'{c}, \emph{Tight contact structures and
  {S}eiberg-{W}itten invariants}, Invent. Math. \textbf{129} (1997), no.~3,
  509--525. \MR{1465333}

\bibitem[LOT18]{LOT_bordered_HF}
Robert Lipshitz, Peter~S. Ozsv\'{a}th, and Dylan~P. Thurston, \emph{Bordered
  {H}eegaard {F}loer homology}, Mem. Amer. Math. Soc. \textbf{254} (2018),
  no.~1216, viii+279. \MR{3827056}

\bibitem[LS22]{lipshitz-sarkar:mixed}
Robert Lipshitz and Sucharit Sarkar, \emph{A mixed invariant of nonorientable
  surfaces in equivariant {K}hovanov homology}, Trans. Amer. Math. Soc.
  \textbf{375} (2022), no.~12, 8807--8849. \MR{4504654}

\bibitem[LZ22]{lewark-zibrowius}
Lukas Lewark and Claudius Zibrowius, \emph{Rasmussen invariants of whitehead
  doubles and other satellites}, arXiv:2208.13612 (2022).

\bibitem[OS04a]{os_genusbounds}
Peter Ozsv\'{a}th and Zolt\'{a}n Szab\'{o}, \emph{Holomorphic disks and genus
  bounds}, Geom. Topol. \textbf{8} (2004), 311--334. \MR{2023281}

\bibitem[OS04b]{os_knotinvts}
\bysame, \emph{Holomorphic disks and knot invariants}, Adv. Math. \textbf{186}
  (2004), no.~1, 58--116. \MR{2065507}

\bibitem[OSS17]{OzStSz_concordanceHoms}
Peter Ozsv\'{a}th, Andr\'{a}s~I. Stipsicz, and Zolt\'{a}n Szab\'{o},
  \emph{Concordance homomorphisms from knot {F}loer homology}, Adv. Math.
  \textbf{315} (2017), 366--426. \MR{3667589}

\bibitem[PW21]{petkova2021twisted}
Ina Petkova and Biji Wong, \emph{Twisted {M}azur pattern satellite knots and
  bordered {F}loer theory}, Michigan Mathematical Journal \textbf{1} (2021),
  no.~1, 1--50.

\bibitem[Ras03]{rasmussen_knotcompl}
Jacob~Andrew Rasmussen, \emph{{F}loer homology and knot complements}, ProQuest
  LLC, Ann Arbor, MI, 2003, Thesis (Ph.D.)--Harvard University. \MR{2704683}

\bibitem[Rot95]{rotman}
Joseph~J. Rotman, \emph{An introduction to the theory of groups}, fourth ed.,
  Graduate Texts in Mathematics, vol. 148, Springer-Verlag, New York, 1995.
  \MR{1307623}

\bibitem[Rub90]{ruberman}
Daniel Ruberman, \emph{Seifert surfaces of knots in s4}, Pacific Journal of
  Mathematics \textbf{145} (1990), 97--116.

\bibitem[Sar15]{sarkar2015moving}
Sucharit Sarkar, \emph{Moving basepoints and the induced automorphisms of link
  {F}loer homology}, Algebr. Geom. Topol. \textbf{15} (2015), no.~5,
  2479--2515. \MR{3426686}

\bibitem[Sch21]{knotjob}
Dirk Sch\"{u}tz, \emph{{K}not{J}ob}, Available at
  \url{http://www.maths.dur.ac.uk/~dma0ds/knotjob.html}, 2021.

\bibitem[Shi71]{shinohara}
Yaichi Shinohara, \emph{Higher dimensional knots in tubes}, Trans. Amer. Math.
  Soc. \textbf{161} (1971), 35--49. \MR{287559}

\bibitem[Sza]{szabo_knot_floer_calc}
Zoltán Szabó, \emph{Knot {F}loer homology calculator},
  \url{https://web.math.princeton.edu/~szabo/HFKcalc.html}.

\bibitem[Tra82]{trace}
Bruce Trace, \emph{On attaching {$3$}-handles to a {$1$}-connected
  {$4$}-manifold}, Pacific J. Math. \textbf{99} (1982), no.~1, 175--181.
  \MR{651494}

\bibitem[Zem16]{Zemke2016QuasistabilizationAB}
Ian Zemke, \emph{Quasi-stabilization and basepoint moving maps in link floer
  homology}, arXiv: Geometric Topology (2016).

\bibitem[Zem18]{zemke_linkcob}
Ian Zemke, \emph{Link cobordisms and functoriality in link {F}loer homology},
  Journal of Topology \textbf{12} (2018), no.~1, 94–220.

\bibitem[Zib22]{zibrowius:2-tangles}
Claudius Zibrowius, \emph{Heegaard floer multicurves of double tangles}, arXiv
  preprint arXiv:2212.08501 (2022).

\end{thebibliography}
\vspace{1cm}
\end{document}